\documentclass[11pt]{amsart}
\usepackage{hyperref}
\usepackage{amsfonts}
\usepackage{amsmath}
\usepackage{amssymb}
\usepackage{amscd}
\usepackage{graphics}
\usepackage{latexsym}
\usepackage{color}
\usepackage{epsfig}
\usepackage{amsopn}
\usepackage{xypic}
\usepackage{mathrsfs}

\usepackage{booktabs}
\usepackage{longtable}
\theoremstyle{plain}
\newtheorem{lemma}{Lemma}[section]
\newtheorem*{theorem*}{Theorem}
\newtheorem*{lemma*}{Lemma}
\newtheorem*{proposition*}{Proposition}
\newtheorem*{conjecture*}{Conjecture}
\newtheorem*{corollary*}{Corollary}
\newtheorem*{problem*}{Problem}
\newtheorem{theorem}[lemma]{Theorem}

\newtheorem{corollary}[lemma]{Corollary}
\newtheorem{proposition}[lemma]{Proposition}
\newtheorem{problem}[lemma]{Problem}
\newtheorem{question}[lemma]{Question}

\theoremstyle{definition}
\newtheorem{definition}[lemma]{Definition}
\newtheorem{example}[lemma]{Example}
\newtheorem{remark}[lemma]{Remark}
\newtheorem{possibility}{Possibility}

\newcommand{\F}[1]{\mathscr{#1}}
\newcommand{\fto}[1]{\stackrel{#1}{\to}}
\newcommand{\OV}[1]{\overline{#1}}
\newcommand{\U}[1]{\underline{#1}}

\newcommand{\Z}{\mathbb{Z}}

\newcommand{\C}{\mathbb{C}}
\newcommand{\Q}{\mathbb{Q}}
\newcommand{\I}{\F{I}}
\newcommand{\OO}{\mathcal{O}}
\newcommand{\te}{\otimes}
\newcommand{\sm}{\setminus}
\newcommand{\id}{\mathrm{id}}
\newcommand{\cI}{\F I}

\renewcommand{\P}{\mathbb{P}}
\newcommand{\PP}{\mathbb{P}}

\DeclareMathOperator{\Bl}{Bl}

\DeclareMathOperator{\Gr}{Gr}
\DeclareMathOperator{\Bs}{Bs}

\DeclareMathOperator{\Mat}{Mat}
\DeclareMathOperator{\Pic}{Pic}

\DeclareMathOperator{\Eff}{Eff}
\DeclareMathOperator{\Ann}{Ann}
\DeclareMathOperator{\rk}{rk}

\DeclareMathOperator{\Mov}{Mov}

\begin{document}

\date{\today}
\author{Jack Huizenga}
\address{Department of Mathematics\\Harvard University, Cambridge, MA 02143}
\email{huizenga@math.harvard.edu}
\thanks{This material is based upon work supported under a National Science Foundation Graduate Research Fellowship}

\title[Steiner bundles and Hilbert schemes]{Restrictions of Steiner bundles and divisors on the Hilbert scheme of points in the plane}

\begin{abstract}
The Hilbert scheme of $n$ points in the projective plane parameterizes degree $n$ zero-dimensional subschemes of the projective plane.  We examine the dual cones of effective divisors and moving curves on the Hilbert scheme.  By studying interpolation, restriction, and stability properties of certain vector bundles on the plane we fully determine these cones for just over three fourths of all values of $n$.   

A \emph{general Steiner bundle} on $\P^N$ is a vector bundle $E$ admitting a resolution of the form $$0\to \OO_{\P^N}(-1)^s\fto M \OO_{\P^N}^{s+r} \to E\to 0,$$ where the map $M$ is general.  We complete the classification of slopes of semistable Steiner bundles on $\P^N$ by showing every admissible slope is realized by a bundle which restricts to a balanced bundle on a rational curve.    The proof involves a basic question about multiplication of polynomials on $\P^1$ which is interesting in its own right.
\end{abstract}

\maketitle

\newcommand{\spacing}[1]{\renewcommand{\baselinestretch}{#1}\large\normalsize}
 \setcounter{tocdepth}{1}
\tableofcontents

\section{Introduction}

An interesting problem in birational geometry is to determine the various birational models of nice moduli or parameter spaces.  Such models often admit interesting geometric interpretations.  As a first step in this problem, it is useful to describe the dual cones of effective divisors and moving curves.

If $X$ is a projective variety, the Hilbert scheme $X^{[n]}$ parameterizes length $n$ zero-dimensional subschemes of $X$.  When $X$ is a smooth curve or surface, $X^{[n]}$ is well-behaved, and is smooth of dimension $n\dim X$.  In these cases, it is a natural compactification of the open symmetric product $(X^n\sm \Delta)/S_n$ parameterizing distinct collections of $n$ points.  If $X$ is a curve, we in fact have $X^{[n]} = X^n/S_n$.  On the other hand, when $X$ is a surface, $X^{[n]}\to X^n/S_n$ is a resolution of the singularities in the symmetric product.

In this paper, we focus on understanding the cones of effective divisors and moving curves on $X^{[n]}$ in the particular case $X = \mathbb{P}^2_{\C}$. Our results are used in \cite{ABCH} to help describe the various birational models of $\P^{2[n]}$.

The Picard group of $\P^{2[n]}$ has rank $2$, and is generated over $\Z$ by classes $H$ and $\Delta/2$, where $H$ is the locus of subschemes meeting a given line and $\Delta$ is the locus of singular subschemes.  The cone $\Eff \P^{2[n]}$ of effective divisors of $\P^{2[n]}$ always has $\Delta$ as one edge, and the other edge is spanned by some divisor of the form $aH-\frac b2\Delta.$  We call $a/b$ the \emph{slope} of such a divisor.  Determining this cone seems to be a fairly subtle problem, and the full cone is only known for certain sporadic values of $n$, such as cases where $n$ is very close to a triangular number. 

When $n=(r+1)(r+2)/2$ is a triangular number, it is easy to describe the effective cone.  In these cases, there is an extremal effective divisor on $\P^{2[n]}$ given as the locus of $n$ points which lie on a curve of degree $r$.  Dually, there is a moving curve class given by allowing $n$ points to move in a linear pencil on a smooth curve of degree $r+1$.  

For a general $n$, we can always write $$n = \frac{r(r+1)}{2} + s \qquad (0\leq s \leq r).$$ A pervasive theme of our results is that the qualitative structure of the effective cone of $\P^{2[n]}$ depends predominately on the value of the ratio $s/r\in [0,1]$, and not on the values of $s$ and $r$ themselves.

\subsection{Interpolation for vector bundles} The concept of interpolation for vector bundles plays a key role in the construction of effective divisors on $\P^{2[n]}$.  Let $X$ be a smooth curve or surface.

 \begin{definition}
 A vector bundle $E$ of rank $r$ on $X$ \emph{satisfies interpolation for $n$ points} if the general $\Gamma \in X^{[n]}$ imposes independent conditions on sections of $E$, i.e. if $$h^0(E\otimes \cI_\Gamma) = h^0(E)-rn.$$  We say $E$ \emph{satisfies unique interpolation for $n$ points} if additionally $h^0(E) = rn$.
 \end{definition}
 
In particular, if $E$ satisfies interpolation for $n$ points then $h^0(E)\geq rn$.  If we let $W\subset H^0(E)$ be a fixed general subspace of dimension $rn$, then we obtain a divisor $D_E(n)$ in $X^{[n]}$ described informally as the locus of schemes of length $n$ which fail to impose independent conditions on sections in $W$.  In case $X = \P^2$ we will compute this divisor's class  as $$[D_E(n)] = c_1(E)H - \frac{r}{2}\Delta.$$ Since we are attempting to compute only the cone $\Eff \P^{2[n]}$, it is worth pointing out that this divisor class is a multiple of the class $$\mu(E)H - \frac{1}{2}\Delta,$$ where $\mu(E) = c_1(E)/\rk(E)$ is the \emph{slope} of $E$.
If we believe that the effective divisors on $\P^{2[n]}$ should come from vector bundles in this fashion, then computation of the effective cone boils down to the following question.

\begin{question}
What is the minimum slope of a vector bundle $E$ on $\P^2$ satisfying interpolation for $n$ points?
\end{question}
Notice that the divisor constructed above for a triangular number $n=(r+1)(r+2)/2$ of points is merely the divisor $D_{\OO_{\P^2}(r)}(n)$.  By studying a more general class of vector bundles called \emph{Steiner bundles}, we will be able to construct many more extremal effective divisors.  

\subsection{Multiplication of polynomials on $\P^1$ and Steiner bundles on $\P^N$}

Consider the following basic problem about polynomial multiplication.  Suppose we have an  $(N+1)$-dimensional subspace $V$ of the space $S_r\subset k[u,v]$ of homogeneous polynomials of degree $r$ in $u,v$.  If $W$ is an $\ell$-dimensional subspace of $S_{s-1}$, think of $W$ as filling up the ``fraction'' $\eta(W)=\ell/s$ of $S_{s-1}$ (noting that $\dim S_{s-1} = s$).  The space $V\cdot W$ spanned by products of elements of $V$ and $W$ lies in $S_{r+s-1}$. 

\begin{question}\label{multProblem}
Let $V\subset S_{r}$ be a general series of dimension $N+1$, and fix a positive integer $s$.  Is it true that for every $W\subset S_{s-1}$ we have $\eta(V\cdot W)\geq \eta(W)$?
\end{question}

In other words, for general $V$, does multiplication of an arbitrary series $W$ by $V$ always increase the fraction of the ambient space that is occupied?

Simple examples show that the answer is not always yes.  For instance, if $r/s>N$ and $V\subset S_r$ is any $(N+1)$-dimensional series, then the multiplication map $$V\te S_{s-1} \to S_{r+s-1}$$ cannot be surjective.  Taking $W=S_{s-1}$, we thus have $\eta(V\cdot W)<\eta(W)=1$.

The answer to this question turns out to be intimately related to properties of certain vector bundles on $\P^{N}$.  A \emph{general Steiner Bundle} $E$ on $\P^{N}$ of rank $r$ is a vector bundle admitting a resolution of the form $$0\to \OO_{\P^{N}}(-1)^s\fto M \OO_{\P^{N}}^{s+r}\to E\to 0,$$ where $M$ is a general matrix of linear forms.  In order for $E$ to be locally free, it is necessary and sufficient that either $s=0$ or $r\geq N$.  These are some of the simplest vector bundles on $\P^{N}$, and much is known about them; we refer the reader to Brambilla \cite{Brambilla2} for an interesting discussion of many of their properties.  Recall that the \emph{slope} $\mu(E)$ of a vector bundle $E$ is given by $c_1(E)/\rk E$, and that $E$ is called \emph{semistable} if every subbundle $F\subset E$ has $\mu(F)\leq \mu(E)$.  See \cite{OSS} for other basic facts about vector bundles that we will use throughout the paper.  Observe that the slope of the bundle $E$ given by the above resolution is $\mu(E) = s/r$.  The next result classifies the slopes of semistable Steiner bundles.

\begin{theorem}\label{stableCriterionThm}
Define a function $\rho_{N}$ by $$\rho_{N}(x) = \frac{1}{N-1+\frac{1}{1+x}}, $$ and put $\phi_N = \lim_{i\to \infty} \rho_N^i(0)$, where $\rho^{i+1} = \rho\circ \rho^i$ and $\rho^0=\id$.  Define a set $\Phi_{N}$ by $$\Phi_{N} = \{\alpha:\alpha > \phi_{N}\} \cup \{\rho_{N}^i(0):i\geq 0\} \subset \Q,$$ The set $\Phi_{N}$ consists of all numbers larger than $\phi_{N}$, together with $0$ and all the convergents in the continued fraction expansion of $\phi_{N}$.   

There exists a semistable Steiner bundle of slope $\mu$ on $\P^N$ if and only if $\mu \in \Phi_N$.  
\end{theorem}

We call the numbers $\rho_N^i(0)$ the \emph{exceptional slopes} of semistable Steiner bundles on $\P^N$.  We note that a large portion of the proof of the theorem follows from earlier work of Brambilla \cite{Brambilla1}, \cite{Brambilla2}.  In particular, Brambilla's work can be seen to imply the nonexistence of semistable Steiner bundles on $\P^N$ with $\mu\notin \Phi_N$, and it also shows the existence of semistable Steiner bundles with slope $\mu$ whenever $\mu$ is exceptional.  We will show that every slope $\mu\in \Phi_N$ can be realized by a semistable Steiner bundle on $\P^N$.

\begin{example} The case of $\P^2$ is the most important for this paper, so we write down the set $\Phi_2$ explicitly as $$\Phi_2 = \{\alpha : \alpha> \varphi^{-1}\} \cup \left\{\frac 01,\frac 12,\frac 35,\frac 8{13},\frac{21}{34},\frac{55}{89},\cdots\right\}\subset \Q \qquad \varphi = \frac{1+\sqrt{5}}{2}.$$ The exceptional slopes are ratios of consecutive Fibonacci numbers, and they converge to the inverse of the golden ratio.
\end{example}

We contrast this classification of the slopes of semistable Steiner bundles with the solution to our original question.

\begin{theorem}\label{multThm}
Let $V\subset S_r$ be a general series of dimension $N+1$, and fix an integer $s$.  Assume $r$ and $s$ are coprime.  Then every series $W\subset S_{s-1}$ satisfies $\eta(V\cdot W) \geq \eta(W)$ if and only if $s/r\in \Phi_{N}$.
\end{theorem}

The coprimality assumption is mainly technical;  we will in fact prove the more interesting reverse direction without this assumption.  The proof of the difficult direction essentially constructs a series $V$ with the required properties under the assumption $s/r\in \Phi_N$.  To prove the easier direction, we show that if there is an $(N+1)$-dimensional series $V$ with the required properties for $s,r$, then there is a semistable Steiner bundle on $\P^N$ with slope $s/r$.

\subsection{A restriction theorem for Steiner bundles} Our primary reason for studying Question \ref{multProblem} is that it gives new insight into Steiner bundles beyond the aforementioned semistability result.  If $E$ is a semistable vector bundle on $\P^N$ and $C$ is a general complete intersection curve of sufficiently high degree, then it is known that $E|_C$ will be semistable; various results to this effect have been given by several authors including Mehta and Ramanathan \cite{MehtaRamanathan} and Flenner \cite{Flenner}.  The general theory does not provide good bounds on how large the degree of $C$ must be, however; furthermore, it also does not usually address what happens for specific types of curves, for instance rational curves.  With the help of  Theorem \ref{multThm}, we are able to give the following result.

\begin{theorem}\label{balancedThm}
Let $E$ be a general Steiner bundle on $\P^N$, given by a resolution $$0\to \OO_{\P^N}(-1)^{ks}\to  \OO_{\P^N}^{k(s+r)}\to E \to 0,$$ and assume $\mu(E)\in \Phi_N$. If $f:\P^1\to \P^N$ is a general degree $r$ map and $k$ is sufficiently large, then $f^\ast E$ has balanced splitting type, i.e. $$f^\ast E \cong \OO_{\P^1}(s)^{kr}.$$  In particular, $E$ is semistable.
\end{theorem}

We believe that the theorem should be true for every $k\geq 1$;  we are able to prove this only when $\mu(E)$ is exceptional, however.  The main idea of the proof is to  show that the property that the pullback is balanced corresponds to some general matrix with entries in an $(N+1)$-dimensional series $V\subset S_r$ giving an isomorphism between two vector spaces of polynomials.  We then look at an incidence correspondence consisting of pairs of matrices and vectors in their kernels, and conclude by a dimension count that the general such matrix has no kernel.  The key estimate in the dimension count is provided by Theorem \ref{multThm}.

\subsection{Interpolation for Steiner bundles on $\PP^2$}

Using our restriction result, we are able to show certain twists and/or duals of Steiner bundles satisfy unique interpolation for $n$ points.  This allows us to determine the effective cone of $\P^{2[n]}$ for just over three quarters of all $n$.

\begin{theorem}\label{hilbThm}
Write $$n = \frac{r(r+1)}{2}+s \qquad (0\leq s \leq r);$$ there is a unique such decomposition.  
\begin{enumerate}
\item If $s/r \in \Phi_2$, then the divisor $$D = (r^2-r+s)H -\frac{r}{2}\Delta$$ spans an edge of the effective cone of $\P^{2[n]}$; if $k$ is sufficiently large then $kD$ is the class of $D_E(n)$, where $E$ is a general bundle with resolution$$0\to \OO_{\P^2}(r-2)^{ks}\to \OO_{\P^2}(r-1)^{k(s+r)} \to E\to 0.$$   These divisors are dual to moving curves $\gamma$ on $\P^{2[n]}$ given by allowing $n$ points to move in a linear pencil on a smooth curve of degree $r$.

\item If $1-\frac{s+1}{r+2}\in \Phi_2$ and $s\geq 1$, then the divisor $$D' = (r^2+r+s-1)H-\frac{r+2}{2}\Delta$$ spans an edge of the effective cone of $\P^{2[n]}$; if $k$ is sufficiently large then $kD'$ is the class of $D_F(n)$, where $F$ is a general bundle with resolution
$$0\to F\to \OO_{\P^2}(r)^{k(2r-s+3)}\to \OO_{\P^2}(r+1)^{k(r-s+1)}\to 0.$$ These divisors are dual to moving curves $\gamma'$ on $\P^{2[n]}$ given by allowing $n$ points to move in a linear pencil on a smooth curve of degree $r+2$.
\end{enumerate}
\end{theorem}

\subsection{Moving curves on $\P^{2[n]}$} Theorem \ref{hilbThm} determines the full effective and moving cones roughly when either $0\leq s/r < 1-\varphi^{-1}\approx 0.382$ or $0.618 \approx \varphi^{-1} < s/r \leq 1$.  The middle region, where $s/r \approx 1/2$, has considerably more complicated behavior, and the full answer here remains open.  Roughly when $2/5 < s/r <1/2$, we will construct a better family of moving curves than those already known.

\begin{theorem}\label{movingSimple}
Write $n = r(r+1)/2 +s$, and suppose $0\leq s < r/2$.  There is a moving curve class $\gamma$ on $\P^{2[n]}$ satisfying 
$$\frac{\gamma\cdot \Delta}{\gamma\cdot H} = \frac{2(2r^2-3r+2s+1)}{2r-1}.$$

\end{theorem}

With $\gamma$ as in the theorem, if $E$ is a general vector bundle with resolution $$0\to \OO_{\P^2}(r-3)^s \to \OO_{\P^2}(r-1)^{2r+s-1}\to E\to 0$$ and $E$ has interpolation for $n$ points, then $\gamma\cdot D_E(n) = 0$, so $\gamma$ is an extremal moving curve and $D_E(n)$ is an extremal divisor.  We suspect such bundles $E$ satisfy interpolation so long as they are semistable, which by \cite{Brambilla2} occurs when either $$\sqrt{2}-1< \frac{s}{r-\frac{1}{2}} < 1/2$$ or $s/(r-\frac{1}{2})$ is a convergent in the continued fraction expansion of $\sqrt 2 -1$.  This moving curve has a bigger slope $\gamma\cdot \Delta/\gamma \cdot H$ than the moving curve given by letting $n$ points move in a linear pencil on a smooth curve of degree $r+2$ so long as $s \geq \frac 15(2r-1).$

The theorem relies on a linkage-type argument making use of the Cayley-Bacharach theorem and a general result on the existence of higher order secant planes to curves in projective space.  Specifically, the theorem depends on a detailed study of the incidence correspondence $$\left\{(\Gamma,\Gamma',\Gamma''):\begin{array}{c}\Gamma\cup\Gamma'\cup\Gamma'' \textrm{ is a reduced complete} \\ \textrm{intersection of two $r$-ics,}\\\textrm{and $(\Gamma\cup \Gamma')\cap L = \emptyset$}\end{array}\right\} \subset \P^{2[n]} \times \P^{2[{r^2-n-(r-1)}]} \times L^{[r-1]},$$ where $L\subset \P^2$ is a line.  

\subsection{Structure of the paper}
We begin by proving our results on Steiner bundles in Sections \ref{multSection}-\ref{slopesSection}.  These sections can be read without any knowledge of the Hilbert scheme.  We then review some basic facts about the Hilbert scheme in Section \ref{hilbReview}.  We prove Theorem \ref{hilbThm} as a consequence of our results on Steiner bundles in Section \ref{InterpSection}.  In Sections \ref{SecantFormulaSection}-\ref{movingSection} we discuss our result on the moving cone.

We conclude the paper in Section \ref{remConeSumm} by summarizing our results and seeing how they fit in with several possible natural conjectures for the remaining cones.  A detailed table describing the known and conjectured effective and moving cones for small $n$ can be found in \cite[Appendix A]{thesis}.  

\subsection*{Acknowledgements} I would like to thank Joe Harris for the guidance that went into preparing this work.  I would also like to thank Izzet Coskun for suggesting this project and for his help.

\section{Multiplication on $\P^1$}\label{multSection}

Our first goal is to prove Theorem \ref{multThm}.  For actually proving the theorem, a renormalization of our notation will be useful.  Recall that we write $S_a = H^0(\OO_{\P^1}(a))$; we choose an affine coordinate  $u$ on $\P^1$, so $S_a$ corresponds to polynomials of degree at most $a$ in $u$.  To avoid trivialities, we will assume $b>a$ and $N \geq 2$ throughout this section.  
If $V\subset S_{b-a}$ is an $N$-dimensional series and $W\subset S_{a-1}$ is a nonempty series, we define the \emph{filling ratio of $W$ with respect to $V$} by $$\mu_V(W) = \frac{\dim (V\cdot W)}{\dim W},$$  where $V\cdot W$ denotes the image of $V\te W\to S_{b-1}$.   In terms of filling ratios, the theorem aims to classify when $\mu_V(W) \geq b/a$ holds for every $W\subset S_{a-1}$ when $V\subset S_{b-a}$ is a general fixed series of dimension $N$.

In the introduction we had $r= b-a$ and $s=a$, so $$\frac{b}{a} = 1+\left(\frac{s}{r}\right)^{-1}.$$ We thus define a set $$\Psi_{N} = 1+\Phi_{N-1}^{-1},$$ defining arithmetic options on sets elementwise, and note that $b/a\in \Psi_N$ if and only if $s/r\in \Phi_{N-1},$ where we interpret division by zero as yielding $\infty$.  The set $\Psi_N$ has a nicer description than $\Phi_{N-1}$ does:  if we put $$\theta(x) = N - x^{-1}$$ and $\psi_N = \lim_{i\to\infty} \theta^i(\infty)$, where we interpret $\theta(\infty)$ as $N$, then it is trivial to verify$$\Psi_N = \{\alpha :1< \alpha < \psi_N\} \cup \{\theta^i(\infty):i\geq 0\} \subset \Q\cup \{\infty\}.$$  We remark that $$\psi_N = \frac{N + \sqrt{N^2-4}}2,$$ so $N-1\leq \psi_N < N$.  Furthermore, every finite element of $\Psi_N$ is no larger than $N$.

Notice that to prove the theorem it suffices to find a single $N$-dimensional  $V$ with the required property.  The next theorem refines one direction of Theorem \ref{multThm}, and its statement will be a bit easier to work with.

\begin{theorem}\label{multThm2}
Suppose $b/a\in \Psi_N$, and let $V\subset S_{b-a}$ be a general series of dimension $N$.  For every nonempty $W\subset S_{a-1}$ we have $\mu_V(W)\geq b/a$.
\end{theorem}

\begin{proof} The proof proceeds in two steps.  First, we will show in Proposition \ref{monomialProp} below that the theorem is true when $1 < b/a\leq N-1$ via a direct argument with monomials.  The theorem is also vacuously true when $a=0$ so that $b/a=\infty$.  Next, if $N-1< b/a\leq N$ and $b/a\in \Psi_N$, put $a' = Na - b$ and $b' = a$.  We will show in Lemma \ref{reductionLemma} that proving the theorem for $a$, $b$, and $N$ can be reduced to proving the theorem for $a'$, $b'$, and $N$.  Notice that the ratio $b'/a'$ satisfies $$\theta\left(\frac{b'}{a'}\right)=N-\frac{a'}{b'} = \frac{b}{a}, \qquad \textrm{so}  \qquad \frac{b'}{a'} = \theta^{-1}\left(\frac ba\right).$$ 
Now look at the function $$\theta^{-1}(x) = \frac{1}{N-x}.$$  We observe that $\theta^{-1}$ has a fixed point at $\psi_N$ (this explains the essential nature of $\psi_N$ to the theorem), and that repeated application of $\theta^{-1}$ will eventually decrease any ratio $b/a$ with $N-1<b/a<\psi_N$ to a ratio $\theta^{-n}(b/a)$ with $1<\theta^{-n}(b/a)\leq N-1$, where the theorem is already known to hold.  On the other hand, if $b/a\in \Psi_N$ and $b/a>\psi_N$, then $b/a = \theta^i(\infty)$ for some $i$, and applying $\theta^{-i}$ reduces us to the trivial case of $b/a=\infty$, completing the proof. 
\end{proof}

On a first reading, it may make sense to skip to the next section at this point, as what follows is both self-contained and one the most technical portions of the paper.   We now proceed to prove the two results cited in the previous proof;   we first show that the theorem holds when $1< b/a \leq N-1$.  All the difficulty of the result occurs already in case $N=3$, so we focus on this case first.

\begin{lemma}\label{ba2}
Suppose $1< b/a\leq 2$, and let $c$ be the remainder upon division of $a$ by $b-a$. The net $$V=\langle 1,u^c,u^{b-a}\rangle$$ satisfies $\mu_V(W)\geq b/a$ for every nonempty $W\subset S_{a-1}$.
\end{lemma}
\begin{proof}
Let $W\subset S_{a-1}$, and consider the space $W'\subset S_{a-1}$ spanned by leading terms (with respect to $u$) of polynomials in $W$.  Clearly $\dim W = \dim W'$.  When we multiply a monomial in $V$ by a monomial in $W'$, we obtain a monomial which is the leading term of an element of $V\cdot W$.  This implies that $\dim(V \cdot W')\leq \dim(V\cdot W)$, and therefore $\mu_V(W')\leq \mu_V(W)$.  Thus to prove the result, we may assume $W$ is spanned by monomials.

We now rephrase the question in terms of sumsets.  Given a set $S\subset \{0,\ldots,a-1\}$, we define the filling ratio of $S$ by $$\mu(S) = \frac{|S+\{0,c,b-a\}|}{|S|},$$ where a sum $S+T$ of two sets of integers denotes $\{s+t:s\in S,t\in T\}$.    We must show $\mu(S) \geq b/a$ for any nonempty $S$.
  
We first reduce to the case where $a,b$ are coprime.  If $k|a$ and $k|b$ then $k|(b-a)$ and $k|c$.  It is easy to see that if the result holds for $a/k$ and $b/k$ then it holds for $a$ and $b$; one can partition $\{0,\ldots,a-1\}$ into the sets $$\{0,k,\ldots,a-k\}\cup \{1,k+1,\ldots,a-k+1\}\cup\cdots \cup \{k-1,2k-1,\ldots, a-1\},$$ and addition of $\{0,c,b-a\}$ respects this decomposition.  

Now assuming $a$ and $b$ are coprime, first suppose that the natural map $\alpha:S\to \Z/(b-a)\Z$ is surjective, so that $S$ contains an integer of each residue class mod $b-a$.  Then $|S+\{0,b-a\}|\geq |S| + b-a$, since $|S+\{0,b-a\}|$ contains a new element in each residue class mod $b-a$.  But $|S|\leq a$, so we conclude  $$\mu(S) = \frac {|S+\{0,c,b-a\}|}{|S|}\geq \frac{|S+\{0,b-a\}|}{|S|}\geq 1+\frac{(b-a)}{|S|}\geq \frac b a.$$

Next assume $\alpha$ is not surjective.  Think of $\Z/(b-a)\Z$ as a graph by joining two residues by an edge whenever they differ by $c$.  Since $c$ is relatively prime to $b-a$, this graph is a connected cycle on $b-a$ vertices.  If the induced subgraph $\alpha(S)$ is not connected, one of its connected components $T\subset \Z/(b-a)\Z$ must satisfy $\mu(\alpha^{-1}(T))\leq \mu(S)$.  Indeed, if $T$ is a component of $\alpha(S)$ and $\OV T$ is the complement of $T$ in $\Z/(b-a)\Z$, then by construction the sets $\alpha^{-1}(T)+\{0,c,b-a\}$ and $\alpha^{-1}(\OV T)+\{0,c,b-a\}$ are disjoint and have union $S+\{0,c,b-a\}$, so $$\mu(S) = \frac{|\alpha^{-1}(T)|}{|S|}\mu(\alpha^{-1}(T)) + \frac{|\alpha^{-1}(\OV T)|}{|S|}\mu(\alpha^{-1}(\OV T)) $$ is the weighted average of $\mu(\alpha^{-1}(T))$ and $\mu(\alpha^{-1}(\OV T))$.  Thus at least one of these numbers is no larger than $\mu(S)$.  Continuing to break up $\OV T \cap \alpha(S)$ into components if necessary, we eventually find a component with the desired property.  We may thus assume that $\alpha(S)$ is connected.

Now that $\alpha(S)$ is connected, it must look like an arithmetic progression with step size $c$:  $$\alpha(S) = \{d,d+c,d+2c,\ldots, d+kc\} \pmod {b-a},$$ where $k$ is between $0$ and $b-a-2$; the above listed elements are all distinct.  We can approximate $$|S+\{0,c,b-a\}| \geq |S|+k+2,$$ since $S+\{0,b-a\}$ contains at least $k+1$ elements not in $S$ (one in each residue class mod $b-a$ in $\alpha(S)$) and $S+c$ has an element whose residue mod $b-a$ has class $d+(k+1)c$, which is not a residue of any element of $S+\{0,b-a\}$.

The last ingredient we need to bound the filling ratio of $S$ is an upper bound on its size.  If $\beta:\{0,\ldots,a-1\}\to \Z/(b-a)\Z$ is the residue map, we can say that $|S|\leq |\beta^{-1}(\alpha(S))|$.  We write $a = (b-a)q+c$ as in the division algorithm.  The fiber of $\beta$ over a residue $e$ in $\Z/(b-a)\Z$ has size $q$ or $q+1$:  it is $q+1$ if $0\leq e<c$, and it is $q$ otherwise.  We must therefore determine how many $h$ of the residues $e$ in $\alpha(S)$ satisfy $0\leq e<c$.

Instead of thinking about residues, think about integers.  Starting at each multiple of $b-a$ we place a ``bucket'' $c$ integers wide, and we are asking how many terms in our arithmetic progression with step size $c$ land in the buckets.  Since the step size of the progression is the same as the bucket width, each bucket can contain at most one term from the progression, and it is impossible to ``skip over'' a bucket.  The arithmetic progression will therefore hit as many buckets as possible if we have $d = c-1$, so that the progression starts at the rightmost edge of a bucket.  The number of buckets hit will equal one more than the number of times the sequence passes a multiple of $b-a$. Therefore $$h \leq 1+\frac{c-1+kc}{b-a}<1+(k+1)\frac{c}{b-a}.$$  We conclude $$|S| \leq |\beta^{-1}(\alpha(S))| = q(k+1)+h <1+(k+1)\left(q+\frac c{b-a}\right).$$

Finally, we finish the proof by observing 
\begin{eqnarray*}
\mu(S) &=& \frac{|S+\{0,c,b-a\}|}{|S|} \geq \frac{|S|+k+2}{|S|}
 = 1+\frac{k+2}{|S|}\\&>& 1+\frac{(b-a)(k+2)}{(b-a)+(k+1)((b-a)q+c)}
= \frac{bk+3b-2a}{ak+b}.
\end{eqnarray*}
But $$\frac{bk+3b-2a}{ak+b} \geq \frac ba,$$ since cross-multiplying shows that it is equivalent to $$(b-a)(2a-b) \geq 0,$$ which is true by assumption.  We conclude $\mu(S) > b/a$, as was to be shown.
\end{proof}

The equivalent result for $N>3$ follows readily from the result for $N=3$, as we will now demonstrate.

\begin{proposition}\label{monomialProp}
Theorem \ref{multThm2} holds when $1< b/a \leq N-1$.
\end{proposition}
\begin{proof}
Write $b-a = qa + r$, choosing the remainder in the range $0< r\leq a$.  Let $V'\subset S_r$ be a net such that for every $W\subset S_{a-1}$ we have $\mu_{V'}(W) \geq (r+a)/a$; this is possible by the lemma since $1<(r+a)/a\leq 2$.  Define $$V = \langle1,u^a,u^{2a},\cdots,u^{(q-1)a}\rangle + u^{qa} V'.$$   Since $$q = \frac ba - 1-\frac ra < \frac ba - 1 \leq N-2$$ we find $$\dim V \leq q + 3 < N+1,$$ and therefore $\dim V \leq N$.  But for $W\subset S_{a-1}$ we have $$V\cdot W \cong W\oplus u^a W\oplus u^{2a}W\oplus\cdots \oplus u^{(q-1)a}W \oplus u^{qa}(V'\cdot W)$$ since the polynomials in $W$ have degree smaller than $a$.  Thus $$\dim(V\cdot W) = q\dim W + \dim(V'\cdot W),$$ and $$\mu_V(W) = q+\mu_{V'}(W) \geq q+1+\frac ra=\frac ba,$$ completing the proof.
\end{proof}

We now complete the second step of the proof of the theorem.

\begin{lemma}\label{reductionLemma}
Put $a'=Na-b$ and $b'=a$, and assume $N-1<b/a \leq N$.  If Theorem \ref{multThm2} holds for $a',b',N$, then it holds for $a,b,N$.
\end{lemma}
\begin{proof}
Consider an inclusion of vector bundles $$0\to \OO_{\P^1}(a'-1)\oplus \OO_{\P^1}(-1)^{N-2} \fto{M} \OO_{\P^1}(a-1)^{N}\to Q\to 0$$
given by a general matrix $M$ of polynomials, and let $Q$ be the cokernel.  Since $N\geq 2$, we find that $Q$ is locally free, hence equals $\OO_{\P^1}(b-1)$ by a Chern class calculation.  We thus have an exact sequence on global sections $$0\to S_{a'-1}\fto\alpha S_{a-1}^{N}\fto\beta S_{b-1}\to 0.$$ Here the map $\beta$ is specified by elements of an at most $N$-dimensional series $V\subset S_{b-a}$; these polynomials are the $(N-1)\times (N-1)$-minors of the matrix $M$.  On the other hand, the map $\alpha:S_{a'-1}\to S_{a-1}^N= S_{b'-1}^N$ is given by independent elements of a general $N$-dimensional series $V'\subset S_{b'-a'}$, so by assumption we may assume $V'$ satisfies the conclusion of the theorem for $a',b',N$.  We claim $V$ satisfies the conclusion of the theorem for $a,b,N$. 

To see this, suppose $W\subset S_{a-1}$ is chosen such that the filling ratio $\mu_V(W)$ is minimal.  If $\mu_V(W) = N$ then we are done, so we may assume $\mu_V(W)<N$, which is to say that $\beta|_{W^N}:W^N\to S_{b-1}$ is not injective.  Write $K = \alpha(S_{a'-1}) = \ker \beta$.  Then $W^N\cap K$ is non-empty.  Let $W'\subset W$ be the subseries spanned by entries of elements of $W^N \cap K$.  Then by construction $$W^N \cap K = (W')^N \cap K.$$ For any series $W\subset S_{a-1}$ we have an exact sequence $$0\to W^N \cap K \to W^N \to V\cdot W\to 0,$$ so $$\dim (V\cdot W) = N\dim W - \dim(W^N\cap K)$$ and $$\mu_V(W) = N - \frac{\dim(W^N\cap K)}{\dim W}.$$  Thus $$\mu_V(W) = N- \frac{\dim(W^N\cap K)}{\dim W} \geq N - \frac{\dim((W')^N\cap K)}{\dim W'} = \mu_V(W'),$$ with equality if and only if $W=W'$.  Since $W$ was chosen with minimal filling ratio, $W=W'$, i.e. $W$ is spanned by the entries of elements of $W^N\cap K$.

Now put $U = \alpha^{-1}(W^N\cap K) \subset S_{a'-1}$.  Clearly $\dim U = \dim (W^N\cap K)$ since $\alpha$ maps $S_{a'-1}$ isomorphically onto $K$.  By the previous paragraph, we see that $V'\cdot U = W$ since $V'\cdot U$ contains all the entries of any element of $W^N\cap K$. 

Finally, since the result holds for $V'$ we have \begin{eqnarray*}\frac{1}{N- \frac ba} =\frac{b'}{a'}& \leq & \mu_{V'}(U) = \frac{\dim(U\cdot V')}{\dim U} = \frac{\dim W}{\dim(W^N\cap K)}\\&=&\frac{\dim W}{N\dim W - \dim(V\cdot W)}=\frac{1}{N-\mu_V(W)},\end{eqnarray*} and we conclude $\mu_V(W) \geq b/a$.
\end{proof}

\section{Matrices with entries in a fixed series}\label{matrixSec}

In this section we prove a result which gives the main link between Steiner bundles and our polynomial multiplication question.

\begin{proposition}\label{netMatrixThm}
Let $V\subset S_{b-a}$ be a general series of dimension $N$, and let $M$ be a general $ak\times bk$ matrix with entries in $V$.  Assume $b/a\in \Psi_N$.  If $k$ is sufficiently large, then the map $$S_{a-1}^{bk}\fto M S_{b-1}^{ak}$$ is an isomorphism.
\end{proposition}
\begin{proof}
We first show that it suffices to consider the case where $a,b$ are coprime.  For say $a=a'd$, $b=b'd$, with $(a',b')=1$.  We can decompose $$S_{a-1} \cong S_{a'-1}^d \qquad S_{b-1} \cong S_{b'-1}^d,$$ where the $i$th factor of each decomposition is spanned by all monomials $u^c$ with $c\equiv i \pmod d$.  If we have a series $V'\subset S_{b'-a'}$ which proves the theorem for $a',b'$, then we can regard it as a series $V\subset S_{b-a}$ by making the change of variables $u\mapsto u^d$.  Then a general matrix $M$ with entries in $V$ will respect the decompositions  $$S_{a-1}^{bk} \cong (S_{a'-1}^{bk})^d \qquad S_{b-1}^{ak}\cong (S_{b'-1}^{ak})^d$$ and give an isomorphism $$S_{a'-1}^{bk} \fto \cong S_{b'-1}^{ak}$$ on each of the $d$ factors separately.  Thus $M$ is an isomorphism.

We now assume $a$ and $b$ are coprime.  Choose $V$ so that the conclusion of Theorem \ref{multThm2} is satisfied.  Observe that $S_{a-1}^{bk}$ and $S_{b-1}^{ak}$ both have dimension $abk$, so to show some $M$ is an isomorphism, it suffices to show it is injective.  Consider the incidence correspondence
$$\xymatrix{
& \ar[dl]_\alpha \Sigma = \{(M,\mathbf{G}):M\mathbf{G}=0\} \ar[dr]^\beta \\
\Mat_{ak\times bk}(V) && \P S_{a-1}^{bk}
}$$ where $\Mat_{ak\times bk}(V)$ denotes the space of $ak\times bk$ matrices with entries in $V$.  We would like to prove that $$\dim \Sigma< \dim \Mat_{ak\times bk}(V)=Nabk^2,$$ since then $\alpha$ is not dominant and the general matrix $M$ gives an isomorphism.  We estimate the dimension of $\Sigma$ by looking at the projection $\beta$.  For $\mathbf{G}\in \P S_{a-1}^{bk}$, we denote by $W_\mathbf{G}$ the subspace of $S_{a-1}$ spanned by the entries of $\mathbf{G}$.  We put $$X_\ell = \{\mathbf{G}:\dim W_\mathbf{G}\leq \ell\}\subset \P S_{a-1}^{bk},$$ and we easily compute $$\dim (X_\ell\sm X_{\ell-1}) = \dim \Gr(\ell,S_{a-1})+bk\ell-1 = \ell(a-\ell)+bk\ell-1.$$   We decompose $$\Sigma = \bigcup_{\ell = 1}^a \beta^{-1}(X_\ell\sm X_{\ell-1}),$$ so we must show each $\beta^{-1}(X_\ell\sm X_{\ell-1})$ has dimension smaller than $Nabk^2$.

To analyze the dimension of $\beta^{-1}(X_\ell\sm X_{\ell-1})$, we must bound the dimension of the fiber over a point $\mathbf{G}\in X_\ell\sm X_{\ell-1}$. If $\mathbf{G}=(g_1,\ldots,g_{bk})$, then a matrix $M$ satisfies $M\mathbf{G} = 0$ exactly when each of its $ak$ rows are in the kernel of \begin{eqnarray*}V^{bk}&\to& S_{a-1}\\(f_1,\ldots,f_{bk})&\mapsto & f_1g_1+\cdots + f_{bk}g_{bk}\end{eqnarray*} The image of this map is $V\cdot W_\mathbf{G}$, so the kernel has dimension $Nbk-\dim (V\cdot W_\mathbf{G})$.  Thus the fiber of $\beta$ over $\mathbf{G}$  has dimension $$\dim \beta^{-1}(\mathbf{G}) = (Nbk-\dim(V\cdot W_\mathbf{G}))ak.$$  Now if $W\in \Gr(\ell,S_{a-1})$ is chosen to minimize $\dim(V\cdot W)$, we estimate $$\dim \beta^{-1}(X_\ell\sm X_{\ell-1}) \leq \ell(a-\ell) +bk\ell -1 + Nabk^2-ak \dim(V\cdot W).$$ We need this quantity to be smaller than $Nabk^2$, which amounts to saying \begin{equation}\label{dimEstimate} \ell(a-\ell) - 1 < k(a\dim (V\cdot W)-b\ell).\end{equation} If $\ell =a$ then this inequality is immediate since $V\cdot S_{a-1}=S_{b-1}$.  Otherwise, if $\ell <a$, we know $\mu_V(W)\geq b/a$.  This inequality is in fact strict, since $$\mu_V(W) = \frac{\dim(V\cdot W)}{\dim W}$$ has denominator smaller than $a$ and $b/a$ is already written in lowest terms. Thus $$a\dim(V\cdot W)-b\ell >0,$$ and for $k$ sufficiently large Inequality (\ref{dimEstimate}) holds.
\end{proof}

\begin{remark}
We believe the conclusion of the proposition holds even if $k=1$.  The argument given here is not refined enough to prove this, however.  To prove the proposition by this general method for $k=1$, it would be necessary to further stratify the $X_\ell\sm X_{\ell-1}$ into loci of the form $$Y_{r,\ell} = \{\mathbf{G}\in X_{\ell}\sm X_{\ell-1}: \dim (V\cdot W_\mathbf{G}) \leq r\}.$$ Theorem \ref{multThm2} shows that $Y_{r,\ell}$ is empty if $b/a\in \Psi_N$ and $r< b\ell/a$.  More generally, we could ask for an upper bound on the dimension of $Y_{r,\ell}$ for all $r,\ell$, and if this estimate is strong enough the result for $k=1$ would follow.

Since this last question seems interesting in its own right, we phrase it in language that does not involve the notation from the proof of Theorem \ref{netMatrixThm}.

\begin{problem}
Let $V\subset S_a$ be a general linear series of dimension $N$.  Estimate the dimension of $$\{W:\dim(V\cdot W)\leq r\} \subset \Gr(\ell,S_b).$$
\end{problem}
\end{remark}

\section{Semistable pullbacks}\label{pullbackSection}

We are now ready to prove our result on the semistability of pullbacks of Steiner bundles to rational curves.  The main observation is that the splitting type of a vector bundle on a rational curve is easy to detect cohomologically.

\begin{theorem}\label{balancedThm2}
Let $E$ be a general Steiner bundle on $\P^N$, given by a resolution $$0\to \OO_{\P^N}(-1)^{ks}\fto M \OO_{\P^N}^{k(s+r)}\to E \to 0,$$ where $s/r\in\Phi_{N}$ and $M$ is given by a general matrix of linear forms. If $f:\P^1\to \P^N$ is a general degree $r$ map and $k$ is sufficiently large, then $$f^\ast E \cong \OO_{\P^1}(s)^{kr}.$$  
\end{theorem}
\begin{proof}
Suppose $f:\P^1\to \P^N$ is given by a general $(N+1)$-dimensional series $V\subset H^0(\OO_{\P^1}(r))$.  The bundle $f^\ast E$ fits into an exact sequence $$0\to \OO_{\P^1}(-r)^{ks} \fto {f^\ast M} \OO_{\P^1}^{k(s+r)} \to f^\ast E\to 0,$$ and the map $f^\ast M$ is given by a general $ k(s+r) \times ks$ matrix with entries in $V$.  Observe that $c_1(f^\ast E) = ksr$, and thus $$f^\ast E \cong \bigoplus_{i=1}^{kr} \OO_{\P^1}(a_i)$$ for some numbers $a_i$ with $\sum a_i = ksr$.  We will have $f^\ast E \cong \OO_{\P^1}(s)^{kr}$ if and only if $$H^0((f^\ast E)^\vee(s-1))=0.$$  Dualizing the above exact sequence and twisting by $\OO_{\P^1}(s-1)$, we get an exact sequence $$0\to (f^\ast E)^\vee(s-1) \to \OO_{\P^1}(s-1)^{k(s+r)} \to\OO_{\P^1}(s+r-1)^{ks}\to 0,$$ so $H^0((f^\ast E)^\vee (s-1)) = 0$ if and only if $$H^0(\OO_{\P^1}(s-1))^{k(s+r)}\to H^0(\OO_{\P^1}(s+r-1))^{ks}$$  is injective.  But $(r+s)/s\in \Psi_{N+1}$ since $s/r\in \Phi_N$, so Proposition \ref{netMatrixThm} completes the proof.
\end{proof}

As a consequence, we obtain the semistability of the above Steiner bundles.

\begin{corollary}\label{stableCor}
For sufficiently large $k$, the bundles of the previous theorem are semistable.  Thus every slope $\mu \in \Phi_N$ is realized by a semistable Steiner bundle.
\end{corollary}
\begin{proof}
In the notation of the theorem, if  $F\subset E$ is a destabilizing subbundle, then $f^\ast F\subset f^\ast E$ is also a destabilizing subbundle, so $E$ is semistable since $f^\ast E$ is.
\end{proof}

This corollary is our contribution to the proof of Theorem \ref{stableCriterionThm} from the introduction; we will complete the proof in the next section.

\section{Slopes of semistable Steiner bundles}\label{slopesSection}

To complete the proof of Theorem \ref{stableCriterionThm}, we must show that if $\mu\notin \Phi_N$ then there is no semistable Steiner bundle of slope $\mu$. While it is not much of a stretch to derive this result from Brambilla \cite{Brambilla2}, the result there is only stated for $\P^2$.  Furthermore, the basic structure of the argument is interesting, and gives insight into Steiner bundles with slope $\mu< \phi_N$.  We therefore sketch the argument, quoting results from Brambilla when necessary.

First of all, we fix $N$ and let $\{a_n\}$ be the sequence defined recursively by \begin{eqnarray*} a_{-1} &=& 0  \\a_0 &=& 1\\ a_{n+1} &=& (N+1)a_n-a_{n-1}.\end{eqnarray*} For $n\geq 0$, we define the \emph{Fibonacci bundle} $F_n$ to be the general Steiner bundle with resolution $$0\to \OO_{\P^N}(-1)^{a_{n-1}} \fto M \OO_{\P^N}^{a_n} \to F_n\to 0.$$  The bundles $F_n$ are exceptional (see \cite{Brambilla1}), so the isomorphism class of $F_n$ is constant as $M$ varies in an open set.  

The following result follows from a trivial induction on $n$.

\begin{lemma}
We have $\mu(F_n) = \rho^n_N(0)$.$\hfill\qed$
\end{lemma}

It is worth recalling that $\rho_N^n(0)$ is an increasing sequence that converges to $\phi_N$.  The main result we will need from Brambilla \cite{Brambilla2} is the following theorem concerning the structure of unstable general Steiner bundles.

\begin{theorem}[Theorem 6.3\cite{Brambilla2}]\label{unstableBundles}
Let $E$ be a general Steiner bundle on $\P^N$, and suppose $$\mu(F_n)\leq \mu (E) < \mu(F_{n+1}).$$ There are uniquely determined integers $k_1$ and $k_2$ such that $$E \cong F_n^{k_1}\oplus F_{n+1}^{k_2}.$$
\end{theorem}

\begin{proof}[Proof of Theorem \ref{stableCriterionThm}]
If $E$ is a general Steiner bundle on $\P^N$ with slope $\mu\notin \Phi_N$, then $0< \mu < \phi_N$.  Thus by the lemma and the theorem $E$ must be a direct sum of two bundles of different slopes, and $E$ is not semistable.
\end{proof}

Notice that since any general Steiner bundle with exceptional slope is a direct sum of copies of a single Fibonacci bundle $F_n$, we can conclude from Theorem \ref{unstableBundles} that Theorem \ref{balancedThm2} holds for all $k\geq 1$ in case the slope is exceptional.

\begin{corollary}\label{k1Cor}
Theorem \ref{balancedThm2} holds for all $k\geq 1$ in case $s/r$ is an exceptional slope.
\end{corollary}

Now that we have finished the classification of semistable slopes of Steiner bundles, it is possible to prove converses to Theorem \ref{multThm2} and Proposition \ref{netMatrixThm}.

\begin{corollary}\label{matrixNetThmConverse}
Let $V\subset S_{b-a}$ be a general $N$-dimensional series, and let $M$ be a general $ak\times bk$ matrix with entries in $V$.  If the map $$S_{a-1}^{bk}\fto M S_{b-1}^{ak}$$ is an isomorphism for some $k\geq 1$, then $b/a\in \Psi_N$.
\end{corollary}
\begin{proof}
By the proofs of Theorem \ref{balancedThm2} and Corollary \ref{stableCor}, the hypotheses imply there is a semistable Steiner bundle on $\P^{N-1}$ with slope $a/(b-a)$.  By Theorem \ref{stableCriterionThm}, $a/(b-a)\in \Phi_{N-1}$ and thus $b/a\in \Psi_N$.
\end{proof}

\begin{corollary}
Suppose $V\subset S_{b-a}$ is an $N$-dimensional series such that $\mu_V(W)\geq b/a$ for every $W\subset S_{a-1}$, where $a$ and $b$ are coprime.  Then $b/a\in \Psi_N$.
\end{corollary}
\begin{proof}
By the proof of Proposition \ref{netMatrixThm}, the general map $$S_{a-1}^{bk}\to S_{b-1}^{ak}$$ given by a matrix with entries in $V$ is an isomorphism for $k$ sufficiently large.  By the previous corollary, $b/a\in \Psi_N$.
\end{proof}

\section{Preliminary facts on the Hilbert scheme of points in the plane}\label{hilbReview}

Here we recall some basic facts about divisor and curve classes on the Hilbert scheme $\P^{2[n]}$ of $n$ points in the plane.  Details for the results in this section can be found in \cite[Chapter 3]{thesis}.

By Fogarty \cite{Fogarty1}, $\P^{2[n]}$ is a smooth projective variety of dimension $2n$.  Since $\mathbb{P}^2$ has irregularity $q=0$, Fogarty \cite{Fogarty2} shows that $$\Pic \P^{2[n]} = \Z H \oplus \Z (\Delta/2),$$ where $H$ is the locus of schemes meeting a fixed line and $\Delta$ is the locus of nonreduced schemes.

Dually, consider the following curves on $\P^{2[n]}$, each parameterized by a $\mathbb{P}^1$.
\begin{itemize}
\item $\alpha$ is the locus where $n-1$ points are fixed and the $n$th point moves on a fixed line.

\item $\beta$ is the locus where $n-2$ points are fixed and a ``spinning tangent direction'' is supported at another fixed point.
\end{itemize}
These classes are a basis for the space $N_1(\P^{2[n]})$ of numerical equivalence classes of curves.  The intersection pairing between divisors and curves is given by the table
$$\begin{array}{c|cc}
&H&\Delta\\
\hline
\alpha & 1 & 0 \\
\beta & 0 & -2
\end{array}$$

Observe that $\alpha$ is a moving curve class, in the sense that the general point of $\P^{2[n]}$ lies on an irreducible curve numerically equivalent to $\alpha$.  Since $\alpha\cdot \Delta = 0$, we conclude that $\Delta$ spans an extremal ray of $\Eff \P^{2[n]}$ and $\alpha$ spans the dual extremal ray of $\Mov \P^{2[n]}$.

Our main task therefore is to determine the other, nontrivial, edge of $\Eff \P^{2[n]}$.  Since $\P^{2[n]}$ is a Mori dream space, this cone is in fact closed, and so we can hope that it is possible to explicitly construct effective divisors spanning this edge.

One other general result on the cones $\Eff \P^{2[n]}$ is easy.  If we fix a point $q\in \P^2$, then unioning a scheme $\Gamma\in \P^{2[n]}$ with $q$ defines a rational map $i_q:\P^{2[n]} \dashrightarrow \P^{2[n+1]}$.   Under this map, we have $i_q^\ast H=H$ and $i_q^\ast \Delta = \Delta$, so we can identify $\Pic \P^{2[n]}$ with $\Pic \P^{2[n+1]}$.  With this identification, we have $\Eff \P^{2[n+1]} \subset \Eff \P^{2[n]}$.

\subsection{The effective divisor corresponding to a bundle with interpolation}\label{divClass}

For every $n$ where we know the full effective cone $\Eff \P^{2[n]}$, it is possible to describe the nontrivial edge of $\Eff \P^{2[n]}$ as a locus where interpolation fails for a vector bundle satisfying interpolation.  
 
 Assume $E$ is a vector bundle of rank $r$ on $\P^2$ that satisfies interpolation for $n$ points.  In particular, we have $h^0(E)\geq rn$.  Let $W\subset H^0(E)$ be a general fixed subspace of dimension $rn$.  A scheme $\Gamma$ which imposes independent conditions on sections of $E$ will impose independent conditions on sections in $W$ if and only if  the subspace $H^0(E\otimes \cI_\Gamma)\subset H^0(E)$ is transverse to $W$.  Thus, informally, we obtain a divisor $D_{E,W}(n)$ described as the locus of schemes which fail to impose independent conditions on sections in $W$.  We observe that the class of $D_{E,W}(n)$ will be independent of the choice of $W$, so we will drop the $W$ when it is either understood or irrelevant to the discussion. 

To put the correct scheme structure on $D_{E,W}(n)$ and compute its class, let $\Xi_n$ be the universal family over $\P^{2[n]}$, with maps as in the diagram $$\xymatrix{\Xi_n \ar[r]^\beta \ar[d]_\alpha & \P^2\\ \P^{2[n]}}$$ The locus of schemes which fail to impose independent conditions on sections in $W$ can be described as the locus where the natural map $$W\otimes \OO_{\P^{2[n]}}\to \alpha_\ast \beta^\ast E =: E^{[n]}$$ of vector bundles of rank $rn$ fails to be an isomorphism.  Consequently, it has codimension at most $1$; since the general $Z$ imposes independent conditions on sections in $W$ it is actually a divisor.  Furthermore, its class (when given the determinantal scheme structure) is just $c_1(E^{[n]})$.  By a simple Grothendieck-Riemann-Roch calculation, we conclude $$[D_E(n)] = c_1(E^{[n]}) = c_1(E)H-\frac{r}{2}\Delta.$$

\subsection{The curve corresponding to a linear pencil on a plane curve}\label{curveClass}

Suppose we are given a smooth plane curve $C\subset \P^2$ of degree $r$, a line bundle $L$ of degree $n$ on $C$, and a linear pencil $\F D \subset \P H^0(L)$.  This pencil induces a map $C\to \P^1$.  Viewing this map as a flat family of degree $n$ zero-schemes parameterized by $\mathbb{P}^1$, we obtain a curve $\gamma$ in $\P^{2[n]}$ parameterized by $\mathbb{P}^1$.  Singular members of $\gamma$ correspond to ramification points of the map $C\to \P^1$.  The Riemann-Hurwitz formula therefore suggests
\begin{eqnarray*}
\gamma\cdot H &=& r\\
\gamma\cdot \Delta &=& 2g(C)-2+2n = r(r-3)+2n.
\end{eqnarray*}  
These formulas are in fact correct, even if $\F D$ has base points or the ramification of the map is complicated.  The necessary tangent space calculations to verify this are straightforward but tedious; we refer the reader to \cite[Section 3.3]{thesis} for a full account.

More generally, suppose $C$ is merely a nodal irreducible plane curve of degree $r$ and geometric genus $g$, with normalization $\tilde C\to C$.  If $L$ is a degree $n$ line bundle on $\tilde C$ and $\F D\subset \P H^0(L)$ is a pencil, we obtain a curve $\tilde \gamma : \P^1\to {\tilde C}^{[n]}$ in the Hilbert scheme of the normalization.  If the general member of $\F D$ is supported away from the preimages of the nodes of $C$, then this induces a curve $\gamma$ in $\P^{2[n]}$ with \begin{eqnarray*}\gamma\cdot H &=& r\\ \gamma\cdot \Delta &\geq& 2(g-1+n).\end{eqnarray*}  Equality holds in the second inequality so long as no member of $\F D$ contains the full preimage of a node, so that no new singular members arise when inducing $\gamma$ from $\tilde \gamma$.

\section{Interpolation for bundles on $\P^2$}\label{InterpSection}

In this section, we prove our strongest result on the effective and moving curves of the Hilbert scheme $\P^{2[n]}$, giving a complete description for about $76\%$  of all values of $n$.  The main task is to show that certain twists and/or duals of general Steiner bundles satisfy interpolation.

\begin{theorem}\label{InterpThm}
Write $$n = \frac{r(r+1)}{2}+s \qquad (s\geq 0),$$ and consider a general vector bundle $E$ given by a resolution $$0\to \OO_{\P^2}(r-2)^{ks}\to \OO_{\P^2}(r-1)^{k(s+r)} \to E\to 0.$$ For sufficiently large $k$, $E$ has (unique) interpolation for $n$ points if and only if $s/r\in \Phi_2$.

Alternately, consider a general vector bundle $F$ given by a resolution $$0\to F\to \OO_{\P^2}(r)^{k(2r-s+3)}\to \OO_{\P^2}(r+1)^{k(r-s+1)}\to 0.$$ For sufficiently large $k$, $F$ has (unique) interpolation for $n$ points if and only if $$1-\frac{s+1}{r+2}\in \Phi_2.$$
\end{theorem}

We will focus primarily on showing Theorem \ref{InterpThm} holds for bundles having the form of $E$ in the theorem, then indicate how analogous results for bundles of the form of $F$ are proved.
The following simple lemma plays a key role in showing semistable twisted Steiner bundles satisfy interpolation. 

\begin{lemma}\label{restrictLemma}
With notation as in Theorem \ref{InterpThm}, if $C\subset \P^2$ is a curve of degree $r$, then the induced map $$H^0(E)\to H^0(E|_C)$$ is an isomorphism, and $H^1(E|_C)=0$.  Furthermore, $h^0(E)=krn$.
\end{lemma}
\begin{proof}
Since $H^1(E)=0$ and there is an exact sequence $$0\to E(-r)\to E\to E|_C\to 0,$$ it suffices to show that $H^0(E(-r)) = H^1(E(-r))=H^2(E(-r))=0$.  This follows immediately from the sequence $$0\to \OO_{\P^2}(-2)^{ks}\to \OO_{\P^2}(-1)^{k(s+r)}\to E(-r)\to 0.$$  The final statement is trivial.
\end{proof}

Thus in order to show a bundle $E$ as above has interpolation, we may take the following approach.  Choose some curve $C\subset \P^2$ of degree $r$, and show there are $n$ points $p_1,\ldots,p_n\in C$ such that $h^0(E|_C(-p_1-\cdots-p_n))=0$.  It then follows that $E$ has no nonzero sections vanishing at $p_1,\ldots, p_n$.  By choosing $C$ to be a general rational curve, we may apply Theorem \ref{balancedThm}.

\begin{proposition}\label{interpProp}
With notation as in Theorem \ref{InterpThm}, if $k$ is sufficiently large and $s/r \in \Phi_2$, then $E$ has interpolation for $n$\ points.
\end{proposition}
\begin{proof}
Since $h^0(E) = krn,$ we must only show that no nonzero sections of $E$ vanish at general points $p_1,\ldots,p_n$.
Let $C\subset \P^2$ be a general rational curve of degree $r$.  By the lemma, $H^0(E)\to H^0(E|_C)$ is an isomorphism.  Since $C$ is general, it has $(r-1)(r-2)/2$ nodes.  We specialize $(r-1)(r-2)/2$ of our $n$ points onto the nodes of $C$, and specialize the remaining $2r+s-1$ points onto smooth points of $C$.  Denote by $D_1$ the divisor of the nodes of $C$ and by $D_2$ the divisor of the smooth points.

Let $f:\P^1\to C$ be the normalization of $C$, and let $\tilde D_1$ and $\tilde D_2$ be the divisors on $\P^1$ lying over $D_1$ and $D_2$, so that $$\deg \tilde D_1 = 2\deg D_1 \qquad \textrm{and} \qquad \deg \tilde D_2 = \deg D_2.$$ Then $$H^0(f^\ast(E|_C)(-\tilde D_1 - \tilde D_2)) \cong H^0(E|_C(-D_1-D_2)).$$ By Theorem \ref{balancedThm}, $$f^\ast(E(-(r-1))|_C) \cong \OO_{\P^1}(s)^{kr},$$ and therefore $$f^\ast(E|_C) \cong \OO_{\P^1}(r^2-r+s)^{kr}.$$ But $$\deg (\tilde D_1 + \tilde D_2) = r^2-r+s+1,$$ so $H^0(f^\ast(E|_C)(-\tilde D_1-\tilde D_2))=H^0(\OO_{\P^1}(-1)^{kr})=0.$ 
\end{proof}

On the other hand, we can show that if $E$ has interpolation then $E$ is semistable.  The key tool is the following result for curves.  

\begin{lemma}\label{interpImpliesSsProp}
Let $E$ be any vector bundle on a smooth curve $C$ with $h^1(E)=0$.  If $E$ has unique interpolation for $n$ points, then $E$ is semistable.
\end{lemma}
\begin{proof}
If $E$ has rank $r$ and has unique interpolation, then $h^0(E) = rn$ and $$h^0(E\te L) = h^1(E\te L) = 0$$ for a general line bundle $L$ of degree $-n$.  Thus $L$ is \emph{cohomologically orthogonal} to $E$, and $E$ is semistable \cite{Faltings}.  (An elementary argument using Riemann-Roch for vector bundles can also be given.)
\end{proof}

\begin{proposition}\label{interpImpSS}
With notation as in Theorem \ref{InterpThm}, if $E$ has interpolation it is semistable, and thus $s/r\in \Phi_2$.
\end{proposition}
\begin{proof}
If $s/r> 1$ we have already seen that $E$ is semistable (regardless of whether it has interpolation), so we assume $s/r\leq 1$.  Let $p_1,\ldots,p_n\in \P^2$ be general points such that $E$ has no nonzero sections vanishing at $p_1,\ldots,p_n$. Since $s \leq r$, there exists a smooth curve $C$ of degree $r$ that contains $p_1,\ldots,p_n$.   By Lemma \ref{restrictLemma}, we have $$h^0(E|_C) = krn, \qquad h^0(E|_C(-p_1-\cdots - p_n)) = 0, \qquad {\textrm{and}}\qquad h^1(E|_C)=0,$$ so $E|_C$ has unique interpolation and Lemma \ref{interpImpliesSsProp} implies $E|_C$ is semistable.  But then $E$ must also be semistable, as a destabilizing subbundle of $E$ would restrict to a destabilizing subbundle of $E|_C$.
\end{proof}

Finally, we address what happens in the case of kernel bundles $F$.

\begin{proposition}
With notation as in Theorem \ref{InterpThm}, if $k$ is sufficiently large then $F$ has interpolation if and only if it is semistable, i.e. if and only if $1-\frac{{s+1}}{r+2}\in \Phi_2$.
\end{proposition}
\begin{proof}
The proof is almost identical to the one given for $E$ if one lets degree $r+2$ curves play the role of the degree $r$ curves in the proof for $E$.  We also apply Theorem \ref{balancedThm} to the dual of $F$ instead of to $F$ itself.  The only nontrivial point is that a priori we could have $h^1(F)>0$, and thus $h^0(F)> k(r+2)n$; however, assuming semistability holds the analogue of Proposition \ref{interpProp} shows that $h^0(F\te \I_{\Gamma})=0$ for a general collection $\Gamma$ of $n$ points, which then forces $h^0(F) = k(r+2)n$.
\end{proof}

With the interpolation result proved, we can finish the proof of Theorem \ref{hilbThm}.

\begin{proof}[Proof of Theorem \ref{hilbThm}]
From the resolutions 
$$\begin{array}{ccccccccc}
0&\to &\OO_{\P^2}(r-2)^{ks} &\to &\OO_{\P^2}(r-1)^{k(s+r)}&\to &E &\to 0\\
0&\to &F &\to & \OO_{\P^2}(r)^{k(2r-s+3)} &\to &\OO_{\P^2}(r+1)^{k(r-s+1)} &\to 0\end{array}
$$
we compute 
\begin{eqnarray*}c_1(E) &=& k(r^2-r+s)\\ c_1(F) &=&  k(r^2+r+s-1).\end{eqnarray*} By Theorem \ref{InterpThm} and the discussion in Section \ref{divClass}, if $k$ is sufficiently large and $s$ and $r$ are such that $E$ (resp. $F$) is semistable, we get corresponding effective divisors $D_E(n)$ (resp. $D_F(n)$) with classes 
$$\begin{array}{rcl}
[D_E(n)] &=& k\left((r^2-r+s)H-\frac{r}{2}\Delta\right) \\

[D_F(n)] &=& k\left((r^2+r+s-1)H-\frac{r+2}{2}\Delta\right). \\
\end{array}$$

Let $\gamma\subset \P^{2[n]}$ be the curve given by letting $n$ points move in a linear pencil on a smooth curve $C$ of degree $r$.  Since $n< h^0(\OO_{\P^2}(r))$, a general collection of $n$ points lies on a smooth curve of degree $r$.  Furthermore, Riemann-Roch asserts this collection moves in a linear pencil on $C$ since $n>g(C)$.  Thus $\gamma$ is a moving curve on $\P^{2[n]}$.  Likewise, allowing $n$ points to move in a linear pencil on a smooth curve $C'$ of degree $r+2$ also gives a moving curve $\gamma'$ on $\P^{2[n]}$ provided $s>1$ so that $n>g(C')=r(r+1)/2$.  By the discussion in Section \ref{curveClass}, we have 
$$
\begin{array}{rclcrcccl}
\gamma\cdot H &=& r &\qquad&\gamma\cdot \Delta &=& r(r-3)+2n &=& 2(r^2-r+s)\\
\gamma'\cdot H &=& r+2 &&\gamma'\cdot \Delta&=& (r+2)(r-1)+2n&=&2(r^2+r+s-1).
\end{array}
$$
We conclude that $\gamma\cdot D_E(n) = \gamma'\cdot D_F(n) =0$, completing the proof.
\end{proof}

\section{Existence of secant planes to curves}\label{SecantFormulaSection}

The principal tool we will use to construct better moving curves on $\P^{2[n]}$ for some remaining $n$ is the existence of higher secant planes to curves in projective space.

\begin{theorem}\label{secantPlanes}
Let $C\subset \P^s$ be a curve of degree $n$ and genus $g$.  Then $C$ has $d$-secant $(d-r-1)$-planes if when we put \begin{eqnarray*}k&=& s+1-d+r\\ \delta &=& n-g-s\end{eqnarray*} we have $$
\delta \geq  0, \qquad
rk  \leq  d,\qquad \textrm{and}\qquad
(r-\delta)k  \leq  g.
$$
\end{theorem}
If we omit the hypothesis $(r-\delta)k\leq g$, this result appeared in \cite[VIII.4, p. 355]{ACGH} as a consequence of the general secant plane formula; with this omission, however, the result is not true.  For instance, without this hypothesis the theorem would imply a degree $4$ elliptic curve in $\P^3$ possesses trisecant lines, which is false.

This result follows easily from the corresponding result for linear series on $C$, whose statement is actually a bit more applicable to our work.  Let $\F D = \P V \subset \P H^0(L)$ be a $g^s_n$ on $C$, and let $V_d^r\subset C^{[d]}$ be the locus of divisors of degree $d$ which impose at most $d-r$ conditions on $\F D$.

\begin{theorem}\label{secantSeries}
Suppose
$$
\delta \geq  0, \qquad
rk  \leq  d,\qquad \textrm{and}\qquad
(r-\delta)k  \leq  g,
$$
where $k$ and $\delta$ are as in Theorem \ref{secantPlanes}.  Then $V_d^r$ is nonempty.  On the other hand, if $$\delta\geq 0,\qquad rk\leq d, \qquad \textrm{and} \qquad (r-\delta)k>g,$$ then either $V_d^r$ is empty or it does not have the expected dimension $d-rk$.
\end{theorem}
\begin{proof}
If $V_d^r$ is empty or has the expected dimension $d-rk\geq 0$, then the class $v_d^r$ of $V_d^r$ in rational cohomology  is computed by the general secant plane formula $$v_d^r = \sum_{1\leq \beta_1<\cdots < \beta_k \leq k+r} \Delta(\beta)^2 \left(\prod_{i=1}^k \mu(r,k,\delta,i,\beta_i)\right)\theta^{\sum(\beta_i-i)}x^{rk-\sum (\beta_i-i)},$$ where 
$\Delta(\beta)$ is the Vandermonde determinant corresponding to $\beta_1,\ldots,\beta_k$, the function $\mu$ is defined by $$\mu(r,k,\delta,i,\beta_i) = {\delta+i-1\choose r+i-\beta_i}\frac{(r+i-\beta_i)!}{(r+k-\beta_i)!(\beta_i-1)!},$$ $\theta$ is the pullback of the theta-divisor on $J(C)$ via the Abel-Jacobi map $C^{[d]}\to J(C)$, and $x$ is the class of the locus of $\Gamma\in C^{[d]}$ containing a fixed point of $C$ \cite[VIII.4, p. 355]{ACGH}.  The binomial coefficient in $\mu$ is defined for arbitrary integers $n$ and $i$ by the convention 
$${n\choose i} = \begin{cases}\displaystyle\frac{n(n-1)\cdots (n-i+1)}{i!} & \textrm{if $i>0$}\\ 1 & \textrm{if $i=0$}\\ 0 & \textrm{if $i<0$.}\end{cases}$$ 

Fix a sequence $1\leq \beta_1<\cdots<\beta_k\leq k+r$ corresponding to a single term $$\Delta(\beta)^2 \left(\prod_{i=1}^k \mu(r,k,\delta,i,\beta_i)\right) \theta^{\sum (\beta_i - i)}x^{rk-\sum(\beta_i-i)}$$ in the sum for $v^r_d$.  Clearly $\Delta(\beta)^2$ is a nonzero positive number.  For each $i$ with $1\leq i\leq k$, we have $\delta+i-1\geq 0$ since $\delta \geq 0$, so we have $${\delta+i-1\choose r+i-\beta_i} \geq 0$$ and thus $\mu(r,k,\delta,i,\beta_i)\geq 0$ for all $i$.  Notice that this binomial coefficient, and hence $\mu(r,k,\delta,i,\beta_i)$, vanishes precisely when $$\delta-1 < r-\beta_i.$$  Thus the product $\prod_i \mu(r,k,\delta,i,\beta_i)$ is positive so long as $\delta-1\geq r-\beta_i$ for all $i$, which occurs whenever $\delta - 1 \geq r-\beta_1$.  If in fact $$\beta_1 \geq r-\delta+1,$$ then we must have $$\beta_i \geq r-\delta+i$$ for all $i$.  

When equality holds for all $i$, so $\beta_i = r-\delta+i$, the corresponding term reduces to a positive number times $$ \theta^{k(r-\delta)}x^{k\delta}.$$ Since $x$ is ample on $C_d$ and $\theta$ is ample on $J(C)$, this cycle is nonempty and effective so long as $k(r-\delta)\leq g$, so that the power of $\theta$ does not exceed $g$.  When this inequality holds, we conclude that the cycle $v^r_d$ is nontrivial, being a finite nonempty sum of positive effective terms, and thus $V_d^r$ is nonempty.  

On the other hand, if $k(r-\delta)>g$, then we have $$\sum_{i=1}^k (\beta_i-i) \geq k(r-\delta)>g$$ for any sequence $\{\beta_i\}$ such that $\prod_i \mu(r,k,\delta,i,\beta_i) \neq 0$.  Since $\theta^{g+1}=0$, the corresponding term is zero, and thus $v^r_d=0$ when $V_d^r$ has the expected dimension.\end{proof}

\section{Some better moving curves on $\P^{2[n]}$}\label{movingSection}

The goal of this section is to construct a family of highly-sloped moving curves on $\P^{2[n]}$ for certain values of $n$ where we have not yet determined the edge of the effective cone.  

\begin{theorem}\label{movingCurves}
Write $n = r(r+1)/2 +s$, and suppose $0\leq s < r/2$.  If $\Gamma\in \P^{2[n]}$ is general, then there is a curve $C\subset \P^2$ of degree $2r-1$ having $$m = r^2-(r-1)-n$$ nodes and no other singularities, such that $\Gamma$ lies on smooth points of $C$ and $\Gamma$ moves in a linear pencil on the normalization $\tilde C$ of $C$.  If $\gamma\subset \P^{2[n]}$ is the corresponding moving curve class, then \begin{eqnarray*}\gamma\cdot H &=& 2r-1\\ \gamma\cdot \Delta &\geq & 2(2r^2-3r+2s+1),\end{eqnarray*} with equality whenever the pencil on $\tilde C$ has no member containing the full preimage of a node of $C$.
\end{theorem}

Modifying $\gamma$ by adding a rational multiple of the moving curve $\alpha$ if necessary, we can produce a moving curve $\gamma'$ with $$\frac{\gamma'\cdot \Delta}{\gamma'\cdot H} = \frac{2(2r^2-3r+2s+1)}{2r-1},$$ so Theorem \ref{movingSimple} follows.  Note that $$m = \frac{1}{2}(r^2-3r-2s+2),$$ so $0\leq m<n$ since $s\geq 0$.

The key ingredient in the proof of the theorem is the study of a particular correspondence.  Fix a line $L\subset \P^2$, and define  
$$\Sigma = \left\{(\Gamma,\Gamma',\Gamma''):\begin{array}{c}\Gamma\cup\Gamma'\cup\Gamma'' \textrm{ is a reduced complete} \\ \textrm{intersection of two $r$-ics,}\\\textrm{and $(\Gamma\cup \Gamma')\cap L = \emptyset$}\end{array}\right\} \subset \P^{2[n]} \times \P^{2[m]} \times L^{[r-1]},$$ noting that $$n+m+(r-1) = r^2.$$ The next proposition summarizes the relevant properties of $\Sigma$ for the proof of the theorem.

\begin{proposition}\label{SigmaProp}
If $n = r(r+1)/2+s$ with $0\leq s<r/2$, then $\Sigma$
\begin{enumerate}
\item is irreducible,
\item dominates $\P^{2[n]}$, and
\item dominates $\P^{2[m]}$.
\end{enumerate}
\end{proposition}

Let us first show that the proposition implies the theorem.  We recall two facts for use in the proof.

\begin{theorem}[Cayley-Bacharach \cite{CayleyBach}]
Let $C_1,C_2\subset \P^2$ be plane curves of degrees $d$, $e$, and suppose that the intersection $\Gamma = C_1\cap C_2$ is zero-dimensional.  Let $\Gamma'$ and $\Gamma''$ be subschemes of $\Gamma$ residual to one another in $\Gamma$, and set $s = d+e-3$.  If $k\leq s$ is a nonnegative integer, then the dimension of the family of curves of degree $k$ containing $\Gamma'$ (modulo those containing all of $\Gamma$) is equal to the failure of $\Gamma''$ to impose independent conditions on curves of complementary degree $s-k$.
\end{theorem}

While we will only need the Cayley-Bacharach theorem in the classical case where $\Gamma$ is reduced, the full concept of residual schemes plays a role in the proof of Proposition \ref{SigmaProp}.  We recall that the subscheme $\Gamma''$ of $\Gamma$ residual to a subscheme $\Gamma'\subset \Gamma$ is the scheme defined by the ideal sheaf $$\F I_{\Gamma''} = \Ann(\F I_{\Gamma'}/\F I_{\Gamma}).$$ For arbitrary schemes this concept is not well-behaved; for instance if $\Gamma'\subset \Gamma$ then the residual to the residual to $\Gamma'$ in $\Gamma$ need not be $\Gamma'$ again.  However, when $\Gamma$ is Gorenstein (which in particular occurs whenever $\Gamma$ is a zero-dimensional complete intersection) everything works nicely.

The other result we will need describes the minimal resolution of the ideal sheaf of a general collection of $n$ points in $\P^2$.

\begin{theorem}[Gaeta \cite{Eisenbud}]
If $n = r(r+1)/2+s$ with $0\leq s \leq r$, then the ideal sheaf $\I_\Gamma$ of a general $\Gamma\in \P^{2[n]}$ admits a resolution $$0\to \OO_{\P^2}(-r-1)^{r-2s} \oplus\OO_{\P^2}(-r-2)^{s}\to  \OO_{\P^2}(-r)^{r-s+1}\to \I_\Gamma\to 0$$ or 
$$0\to \OO_{\P^2}(-r-2)^s \to \OO_{\P^2}(-r)^{r-s+1}\oplus \OO_{\P^2}(-r-1)^{2s-r} \to \I_\Gamma\to 0,$$
depending on whether $s\leq r/2$ or $s\geq r/2$.  In either case, the homogeneous ideal of $\Gamma$ is generated by $r$-ics and $(r+1)$-ics.
\end{theorem}

\begin{proof}[Proof of Theorem \ref{movingCurves}]
Since $\Sigma$ dominates $\P^{2[n]}$, for a general $\Gamma\in \P^{2[n]}$ we can find a triple $(\Gamma,\Gamma',\Gamma'')\in \Sigma$.  Since $\Gamma\cup \Gamma'$ is linked in a complete intersection of $r$-ics to the collinear collection $\Gamma''$, by Cayley-Bacharach the collection $\Gamma\cup \Gamma'$ fails to impose independent conditions on curves of degree $2r-4$.  Then by Riemann-Roch, if there is a curve $C$ of degree $2r-1$ passing through $\Gamma$, nodal at each point of $\Gamma'$, and having no further singularities, then $\Gamma$ moves in a linear pencil on the normalization $\tilde C$ of $C$.  Riemann-Hurwitz says that this linear pencil on $\tilde C$ has $$2g(\tilde C)-2-n(2g(\P^1)-2) = 2(2r^2-3r+2s+1)$$ singular members (with multiplicity).  If $\gamma\subset\P^{2[n]}$ is the induced curve in $\P^{2[n]}$, then by the discussion in Section \ref{curveClass} \begin{eqnarray*}
\gamma \cdot H &=& 2r-1\\
\gamma\cdot \Delta &\geq & 2(2r^2-3r+2s+1),
\end{eqnarray*}
with equality whenever no additional points of $\gamma\cap \Delta$ arise when the pencil descends from $\tilde C$ to $C$.  We must therefore show that such a curve $C$ exists.

Consider the blowup $\Bl_{\Gamma\cup \Gamma'\cup \Gamma''} \P^2$, and denote by $E,F,G$ the sums of exceptional divisors corresponding to $\Gamma,\Gamma',$ and $\Gamma''.$  By construction, the series $|rH - E - F - G|$ is nonempty and base point free.  If we show that $|(r-1)H - F|$ is nonempty and base point free, then $|(2r-1)H-E-2F-G|$ is base-point free and its general member will be smooth by Bertini.  Furthermore, the general curve of degree $r$ vanishing along $\Gamma\cup \Gamma'\cup\Gamma''$ has general tangent directions at points of $\Gamma'$ (since $\Gamma\cup\Gamma'\cup\Gamma''$ is a transverse complete intersection of $r$-ics), so the general member of $|(2r-1)H-E-2F-G|$ meets each of the components of $F$ in a distinct pair of points.  Thus the general member of this series corresponds to a plane curve of degree $2r-1$  passing through $\Gamma$, nodal at each point of $\Gamma'$, and having no further singularities.

To finish the proof, we must therefore show $|(r-1)H-F|$ is nonempty and base-point free.  Nonemptiness is obvious, as $$\dim |(r-1)H-F| \geq h^0(\OO_{\P^2}(r-1))-1 - m = 2r+s-2,$$ with equality whenever $\Gamma'$ imposes independent conditions on curves of degree $r-1$.  Thus we concentrate on base-point freeness. 

We claim that if $X\subset \P^{2[m]}$ is any proper subvariety, then for general $\Gamma\in \P^{2[n]}$ we may find a triple $(\Gamma,\Gamma',\Gamma'')\in \Sigma$ such that $\Gamma'\notin X$.  Indeed, if $\beta: \Sigma\to \P^{2[m]}$ is the projection, then since it is dominant we see that $U = \beta^{-1}(\P^{2[m]}\sm X)\subset \Sigma$ is a nonempty open subset.  Since $\Sigma$ is irreducible, $U$ is a dense open subset.  But since $\Sigma$ dominates $\P^{2[n]}$, so does $U$, and the claim follows.

In particular, if we take $X\subset \P^{2[m]}$ to be the locus of $\Gamma'$ such that $|(r-1)H-F|$ is not base-point free, it suffices to show that $X$ is a proper subvariety of the Hilbert scheme.  For this it suffices to know that the general $\Gamma'\in \P^{2[m]}$ is cut out scheme-theoretically by $(r-1)$-ics, or that its ideal has a set of generators with degrees no more than $r-1$. But since $m < {r\choose 2},$ Gaeta's theorem implies the ideal of $\Gamma'$ is generated by polynomials of degree at most $r-1$.  
\end{proof}

The most important aspect of $\Sigma$ for the previous proof is the fact that it dominates $\P^{2[n]}$, as this is the condition that ensures we can find a potential location $\Gamma'$ for the nodes of $C$.  Unfortunately the full proof that $\Sigma$ dominates $\P^{2[n]}$ is rather technical, even though the basic idea is simple.  The next lemma contains the key insight of the proof, and also explains the occurrence of the condition $s<r/2$.

\begin{lemma}\label{Gamma''Lemma}
If $\Gamma\in \P^{2[n]}$ is general and $0\leq s < r/2$, then there is some $\Gamma''\in L^{[r-1]}$ such that $\Gamma\cup \Gamma''$ lies on a pencil of $r$-ics.
\end{lemma}
\begin{proof}
Let $V = H^0(\I_{\Gamma}(r))|_L$, and consider the linear series $\F E = \P V$ on $L$.  Every member of $H^0(\I_{\Gamma}(r))$ is irreducible since $\Gamma$ lies on no $(r-1)$-ic, so the restriction map $H^0(\I_{\Gamma}(r))\to V \subset H^0(\OO_L(r))$ is injective.  Thus $$\dim \F E = {r+2\choose 2} - 1 -n = r-s,$$ so $\F E$ is a $g_r^{r-s}$ on the rational curve $L$.  We must show that there is some divisor $\Gamma''$ on $L$ of degree $r-1$ such that $\F E(-\Gamma'')$ has dimension at least $1$; if we can do this then a lift of a pencil in $\F E(-\Gamma'')$ to a pencil in $\P H^0(\I_{\Gamma}(r))$ will be a pencil vanishing on $\Gamma\cup \Gamma''$.

Thus we wish to show that the locus $V_{r-1}^s$ of divisors $\Gamma''\subset L$ of degree $r-1$ which fail to impose at least $r-s$ conditions on $\F E$ is nonempty.  We denote by $\U s$, $\U n$, $\U r$, $\U d$, $\U g$, $\U k$, $\U \delta$ the variables from Section \ref{SecantFormulaSection}, apologizing for the conflicts with our current notation.  Then we have $$\begin{array}{rclcrclcrcl} \U s &=& r-s & \quad& \U n &=& r & \quad & \U r &=& s\\
\U d &=& r-1 && \U g &=& 0  \\ \U k &=& \U s+1 - \U d+ \U r =2 && \U\delta &=& \U n-\U g -\U s = s.\end{array}$$ The three inequalities $$
\U \delta \geq  0, \qquad
\U r \, \U k  \leq  \U d,\qquad \textrm{and}\qquad
(\U r-\U \delta)\U k  \leq \U g
$$
are all satisfied since $0\leq s < r/2$, so in fact $V_{r-1}^s$ is nonempty by Theorem \ref{secantSeries}.
\end{proof}

Now that we have the lemma, assume $\Gamma\in \P^{2[n]}$ is general, and find some $\Gamma''\in L^{[r-1]}$ such that $\Gamma\cup\Gamma''$ lies on a pencil of $r$-ics.  Then we can let $\Gamma'$ be the scheme residual to $\Gamma\cup \Gamma''$ in the base locus $\Bs \F D$.  Assuming that $\Bs \F D$ is reduced and meets $L$ exactly in $\Gamma''$, the triple $(\Gamma,\Gamma',\Gamma'')$ lies in $\Sigma$, and we are done.  Justifying this assumption requires substantial effort, however.

\begin{lemma}\label{SigmaHnLemma}
If $0\leq s<r/2$, then $\Sigma$ dominates $\P^{2[n]}$.
\end{lemma}
\begin{proof}
We must introduce a couple auxiliary correspondences.
First, let $$\F X = \{(\Gamma'',\F D): \F D \in \Gr(2,H^0(\I_{\Gamma''\subset \P^2}(r)))\} \subset L^{[r-1]}\times \Gr(2,H^0(\OO_{\P^2}(r))).$$ That is, $\F X$ is the Grassmannian bundle over $L^{[r-1]}$ corresponding to $2$-planes in the vector bundle on $L^{[r-1]}$ whose fiber over a point $\Gamma''\in L^{[r-1]}$ is $H^0(\I_{\Gamma''\subset \P^2}(r))$.  Since every $\Gamma''\in L^{[r-1]}$ imposes $r-1$ independent conditions on $r$-ics, Grauert's theorem \cite[III.12.9]{Hartshorne} implies this vector bundle can be constructed as a pushforward.  Clearly $\F X$ is projective and irreducible.  We must also consider the dense open subsets in $\F X$ given by \begin{eqnarray*}\F X_1 &=& \{(\Gamma'',\F D):\textrm{$\Bs \F D$ is zero-dimensional}\}\\ \F X_2 &=& \left\{(\Gamma'',\F D):\textrm{$\Bs \F D$ is zero-dimensional, reduced, and $\Bs \F D \cap L = \Gamma''$}\right\}.\end{eqnarray*} It is easy to see that these are in fact dense open subsets of $\F X$.

We now introduce the correspondence $$ \F Y = \left\{(\Gamma,\Gamma''',(\Gamma'',\F D)):\begin{array}{c} \Gamma \textrm{ reduced}\\ \Gamma\cap L = \emptyset \\ \Gamma \textrm{ lies on no $(r-1)$-ic}\\ \Gamma\subset \Gamma''' \\ \Gamma''' \textrm{ is residual to $\Gamma''$ in $\Bs \F D$}  \end{array}\right\}\subset \P^{2[n]} \times \P^{2[n+m]} \times \F X_1,$$ noting that $\F D$ automatically has zero-dimensional base locus since $\Gamma\subset \Bs \F D$ and $\Gamma$ lies on no $(r-1)$-ic, so that every member of $\F D$ is irreducible.

Lemma \ref{Gamma''Lemma} implies that $\F Y$ dominates $\P^{2[n]}$.  We saw that for general $\Gamma\in \P^{2[n]}$ there is some $\Gamma''\in L^{[r-1]}$ and a pencil $\F D$ such that $\Gamma\cup \Gamma'' \subset \Bs \F D$.  Now if we take $\Gamma'''$ to be residual to $\Gamma''$ in $\Bs \F D$, the point $(\Gamma,\Gamma''',(\Gamma'',\F D))$ lies in $\F Y$.

Observe that we have a map $$\phi:\F X_1\to \P^{2[n+m]}\times \F X_1$$ given by sending a point $(\Gamma'',\F D)$ to $(\Gamma''',(\Gamma'',\F D))$, where $\Gamma'''$ is the scheme residual to $\Gamma''$ in $\Bs \F D$.  If $$\alpha:\F Y\to\P^{2[n+m]}\times \F X_1$$ is the projection, then we see that $$\alpha^{-1}(\phi(\F X_1))=\F Y.$$  We claim that $\alpha^{-1}(\phi(\F X_2))$ contains a dense open subset of $\F Y$; if this is true then $\alpha^{-1}(\phi(\F X_2))$ will dominate $\P^{2[n]}$, which implies $\Sigma$ dominates $\P^{2[n]}$.  

Let $(\Gamma,\Gamma''',(\Gamma'',\F D))$ be any point in $\F Y$, and choose a general curve $\gamma(t) = (\Gamma''_t,\F D_t)$ in $\F X$ with $\gamma(0) = (\Gamma'',\F D)$, parameterized by a disc $\Delta$.  For small nonzero $t$, we have $\gamma(t) \in \F X_2$, so that the scheme $\Gamma'''_t$ residual to $\Gamma''_t$ in $\Bs \F D_t$ is a reduced collection of points disjoint from $\Gamma''_t$.  After a base change if necessary, we may assume that there are $n+m$ arcs $p_i(t)$ in $\P^2$ parameterized by $\Delta$ such that $\Gamma'''_t = \{p_1(t),\ldots,p_{n+m}(t)\}$ for nonzero $t$.  Without loss of generality we may assume $\Gamma =\{p_1(0),\ldots,p_n(0)\}$; the assumption that $\Gamma$ is reduced is crucial here.  Then if we write $\Gamma_t = \{p_1(t),\ldots,p_n(t)\}$, we obtain an arc \begin{eqnarray*}\tilde \gamma: \Delta& \to& \F Y\\ t&\mapsto &(\Gamma_t,\Gamma'''_t,(\Gamma''_t,\F D_t))\end{eqnarray*} 
such that $\tilde\gamma (0) = (\Gamma,\Gamma''',(\Gamma'',\F D))$ and $\tilde \gamma(t)$ lies in $\alpha^{-1}(\phi(\F X_2))$ for small nonzero $t$.  Thus $\alpha^{-1}(\phi(\F X_2))$ is dense in $\F Y$, completing the proof.
\end{proof}

The proof that $\Sigma$ dominates $\P^{2[m]}$ is much easier.

\begin{lemma}
If $0\leq s< r/2$, then $\Sigma$ dominates $\P^{2[m]}$.
\end{lemma}
\begin{proof}
Note that $m<n$. If $\Gamma'\in \P^{2[m]}$ is general, choose a general $\Gamma_0 \in \P^{2[n-m]}$.  Then $\Gamma'\cup \Gamma_0$ is general in $\P^{2[n]}$, so since $\Sigma$ dominates $\P^{2[n]}$ we can find some triple $(\Gamma'\cup \Gamma_0, \Gamma_1,\Gamma'')\in \Sigma$. Clearly then $(\Gamma_0\cup \Gamma_1,\Gamma',\Gamma'')\in \Sigma$, so $\Sigma$ dominates $\P^{2[m]}$.
\end{proof}

We conclude the proof of the proposition by showing $\Sigma$ is irreducible.

\begin{lemma}
If $0\leq s < r/2$, then $\Sigma$ is irreducible.  In fact, if we have $$n = r(r+1)/2+s$$ and only assume $s\geq 0$, then $\Sigma$ is irreducible so long as $m\geq 0$, so that the definition of $\Sigma$ makes sense.
\end{lemma}
\begin{proof}Let us put $$\Xi = \{(p,(\Gamma'',\F D)):p\in \Bs \F D\} \subset (\P^2\sm L)\times \F X_2,$$ and observe that $\Xi$ is an $(n+m)$-sheeted covering space of $\F X_2$ (see the proof of Lemma \ref{SigmaHnLemma} for the definition of $\F X_2$).  Denote by $$\Xi(k) = \{(p_1,\ldots,p_k,(\Gamma'',\F D)):\textrm{$p_i\in \Bs \F D$ distinct}\}\subset ((\P^2\sm L)^k\sm \Delta)\times \F X_2.$$ If we show $\Xi(n+m)$ is irreducible, then $\Sigma$ is irreducible since it is the image of $\Xi(n+m)$ under the map $$(p_1,\ldots,p_{m+n},(\Gamma'',\F D)) \mapsto  (\{p_1,\ldots,p_n\},\{p_{n+1},\ldots,p_{n+m}\},\Gamma'').$$ As $\Xi(n+m)$ is an \'etale $(m+n)!$-sheeted cover of the irreducible variety $\F X_2$, it suffices to show that $\Xi(n+m)$ is connected.  Equivalently, we can show that the monodromy group of $\Xi \to \F X_2$ is the full symmetric group.  To do this, we show the monodromy acts doubly transitively on a fiber  and that it contains a simple transposition.
 
To see the monodromy acts doubly transitively, it suffices to show that $\Xi(2)$ is connected.  Consider the slightly enlarged correspondence $$\tilde \Xi(2)=\{(p_1,p_2,( \Gamma'',\F D)):\textrm{$p_i\in \Bs \F D$ distinct}\} \subset ((\P^2\sm L)^2\sm \Delta) \times \F X,$$ observing in particular that we have allowed $\F D$ to be an arbitrary pencil containing $\Gamma''$.  Then the fiber of $\tilde \Xi(2)$ over the triple $(p_1,p_2,\Gamma'')$ consists of the Grassmannian $$\Gr(2,H^0(\I_{\{p_1,p_2\}\cup \Gamma''}(r))).$$  Observing that the dimension $h^0(\I_{\{p_1,p_2\}\cup \Gamma''}(r))$ is independent of $p_1,$ $p_2,$ and $\Gamma'',$ we conclude by Grauert's theorem that $\tilde \Xi(2)$ is actually the Grassmannian bundle associated to a vector bundle on the irreducible variety $(\P^2\sm L)^2 \sm \Delta$.  Thus $\tilde \Xi(2)$ is irreducible, and the dense open subset $\Xi(2)$ is connected.

To find a simple transposition, consider the following family of pencils of $r$-ics, corresponding to a real loop in $\F X_2$.  Fix $r-1$ points $\Gamma''$ on $L$, let $C_1$ be a general curve of degree $r$ containing $\Gamma''$, and let $C_2$ be a general curve of degree $r-1$ containing $\Gamma''.$ Pick a general tangent line to $C_1$, and choose affine coordinates on $\P^2$ so that this line is given by $y=0$ and tangent to $C_1$ at the origin.  Fix a small $\varepsilon$, and denote by $L_t$ the line given by $y=\varepsilon e^{2\pi i t}$. Put $\F D_t = \langle C_1,C_2\cup L_t\rangle.$ As $t$ goes from $0$ to $1$, this gives a loop in $\F X_2$ based at $\langle C_1,C_2\cup L_0\rangle$.  The local equation of $C_1$ near the origin is $y=x^2$, so we see that the induced element of the monodromy group exchanges the pair of points of $L_0 \cap C_1$ near the origin with one another while leaving all other base points of $\F D_0$ fixed.  Thus the monodromy group contains a simple transposition, and is the full symmetric group.
\end{proof}

\section{Summary of results and conjectures}\label{remConeSumm}

Let us briefly summarize an approximate picture of the current results and conjectures on the effective and moving cones of $\P^{2[n]}$.  The actual results are slightly more messy, and are made precise in the indicated places.  In particular, we neglect to mention the sporadic cases corresponding to $s/r$ being a continued fraction expansion of some irrational number, and focus solely on continuous behavior. 

\begin{figure}\label{slopesGraph}
\begin{center}
\input{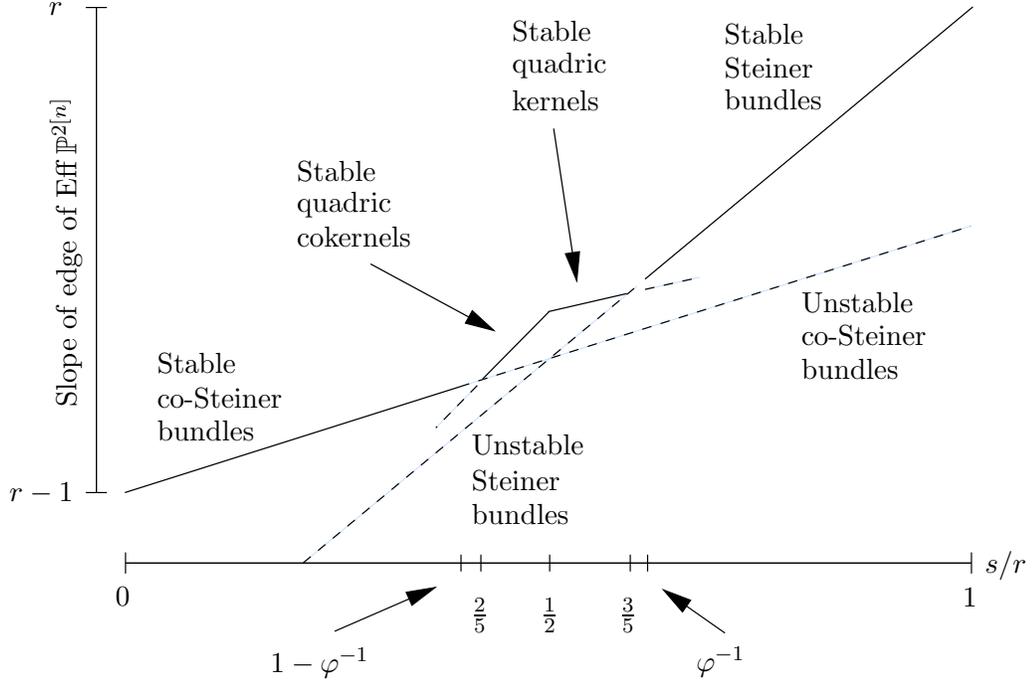}
\end{center}
\caption{Conjectural schematic picture of the slope of the nontrivial edge of the effective cone of $\P^{2[n]}$, where $n = r(r+1)/2+s$, $r$ is fixed, and $s$ ranges from $0$ to $r$.  The image is distorted to emphasize the relative slopes between the lines; for large $r$ these slopes are all very similar.}
\end{figure}

Recall that $$n = \frac{r(r+1)}{2} + s \qquad (0\leq s\leq r).$$  The structure of the effective cone of $\P^{2[n]}$ predominately depends on the value of the ratio $s/r$, with only some slight inaccuracy; in the asymptotic picture as $r$ becomes large, this inaccuracy vanishes.  Figure \ref{slopesGraph} summarizes the following discussion.

\textbf{Case 1.} $0<s/r \lesssim 1-\varphi^{-1} \approx 0.382$. By Theorem \ref{hilbThm},  the nontrivial extremal ray of the effective cone is spanned by $$\mu(E)H - \frac 12 \Delta = \frac{r^2+r+s-1}{r+2} H-\frac{1}{2}\Delta,$$ where $E$ is a vector bundle with resolution
$$0\to E \to \OO_{\P^2}(r)^{2r-s+3}\to \OO_{\P^2}(r+1)^{r-s+1}\to 0,$$ so the extremal ray is given by stable \emph{Steiner kernel bundles.}  Dually, the moving curve is given by letting $n$ points move in a linear pencil on a smooth curve of degree $r+2$.

\textbf{Case 2.} $0.414\approx \sqrt{2} - 1 \lesssim s/r < 1/2$.  There is a moving curve given by letting $n$ points move in a linear pencil on a curve of degree $2r-1$ having $$m = r^2-(r-1)-n$$ nodes and no further singularities---this construction works when $0<s/r<1/2$, but provides a moving curve of slope higher than the moving curve from Case 1 roughly when $s/r >2/5$ (Theorem \ref{movingCurves}).  We conjecture this moving curve is extremal, dual to the divisor given as the locus where interpolation fails for a general vector bundle with resolution $$0\to \OO_{\P^2}(r-3)^s \to \OO_{\P^2}(r-1)^{2r+s-1}\to E \to 0,$$ a stable \emph{quadric cokernel bundle}; interpolation is not known for these bundles, but computer evidence verifies it holds for small $n$.  The extremal ray would be spanned by $$\frac{{2r^2-3r+2s+1}}{2r-1}H-\frac 12 \Delta.$$

\textbf{Case 3.} $1/2 < s/r \lesssim 2-\sqrt{2} \approx 0.586$.  Dually to the previous case, we predict the edge of the effective cone corresponds to divisors coming from stable \emph{quadric kernel bundles}, with resolution of the form $$0\to E \to \OO_{\P^2}(r)^{3r-s+6}\to \OO_{\P^2}(r+2)^{r-s+1}\to 0.$$ The moving curve should be given by allowing $n$ points to move in a linear pencil on a curve of degree $2r+5$ having $$m = (r+3)^2 - (r+2)-n$$ nodes and no other singularities (the existence of this moving curve class is also not known).  The extremal ray is spanned by $$\frac{2r^2+3r+2s-2}{2r+5} H - \frac{1}{2}\Delta.$$ For additional discussion of Cases 2 and 3, see \cite[Section 5.4]{thesis}

\textbf{Case 4.} $0.618 \approx \varphi^{-1}< s/r \leq 1$.  Here Theorem \ref{hilbThm} asserts the edge is spanned by the class $$\frac{r^2-r+s}{r}H-\frac{1}{2}\Delta,$$ corresponding to a stable \emph{Steiner cokernel bundle} with resolution $$0 \to \OO_{\P^2}(r-2)^{s}\to \OO_{\P^2}(r-1)^{r+s} \to E\to 0.$$  Allowing $n$ points to move in a linear pencil on a smooth curve of degree $r$ yields the moving curve.

\begin{remark}\label{possRemark}  While the previous theorems and conjectures address the vast majority of all $n$, approximately $6.4\%$ of all cases are completely open.  There are a few natural guesses as to the slope of the effective cone for the remaining $n$, each generalizing the current data.  For simplicity, let us focus on the case $1/2<s/r\leq 1$; the other cases have a dual picture.
\end{remark}

\begin{possibility}\label{easyPoss} Given an $n$ with $1/2<s/r\leq 1$, we have constructed two different moving curve classes, as in Cases 3 and 4 above.  Perhaps one of these two moving curve classes is always extremal.\end{possibility}

If this is the case, then there exist $n$ such that the dual extremal effective divisors  do not come from vector bundles satisfying unique interpolation.  Specifically, when $3/5< s/r < \varphi^{-1}$ and $s/r\notin \Phi_{2}$, this possibility would predict the edge of the effective cone is spanned by the class $$\frac{r^2-r+s}{r}H-\frac{1}{2}\Delta.$$  However, one can show that if $E$ is a vector bundle satisfying unique interpolation for $n$ points and $\mu(E) = (r^2-r+s)/r$ then $E$ must admit a resolution of the form $$0\to \OO_{\P^2}(r-2)^{ks} \to \OO_{\P^2}(r-1)^{k(s+r)}\to E\to 0$$ for some $k\in \Q$, and thus $s/r\in \Phi_2$.  See \cite[Section 5.5]{thesis} for details.

\begin{possibility}\label{stablePoss}
We may need to allow more general vector bundles.  Noting that in the solved and conjectured cases stable bundles have played an important role, we can choose a general stable vector bundle $E$ having minimal slope among vector bundles with $\chi(E) = n\cdot \rk E$.  If such a bundle has interpolation, it may yield a divisor spanning the edge of the effective cone.
\end{possibility}

Drezet and Le Potier's classification in \cite{DLP} of the possible numerical invariants of stable vector bundles allows one to determine the minimal slope in the preceding construction, and this slope agrees with the conjectured slope in the Steiner and quadric cases.  A difficulty with this possibility is that the ranks of the minimal vector bundles become astronomical (in fact, are unbounded) when $s/r$ is around $3/5$, which means the dual moving curves would correspond to points moving on plane curves of very large degree.

\begin{possibility}\label{homoPoss}
Finally, instead of allowing arbitrary stable vector bundles, perhaps we should only allow stable vector bundles where the minimal resolution of either the bundle or its dual is homogeneous, in the sense that the matrix consists entirely of forms of the same degree.
\end{possibility}

We feel that this option seems somewhat less natural, but it has the upside of allowing divisors to come from vector bundles while bounding the complexity of the situation around $s/r = 3/5$.

\begin{remark}
The first case not handled by Steiner or quadric bundles is $$n=142 = \frac{16\cdot 17}{2}+6$$  Already in this case, all three possibilities disagree.  It is easy to check that Possibilities 1, 2, and 3 predict edges spanned by $$\frac{277}{18}H-\frac 12 \Delta, \qquad \frac{1185}{77} H-\frac{1}{2}\Delta, \qquad\textrm{and} \qquad \frac{77}{5}H-\frac 12 \Delta,$$ respectively.  Note in particular that if Possibility 2 is accurate, the dual moving curve involves 142 points moving on a curve of degree divisible by 77.
\end{remark}

\bibliographystyle{plain}

\end{document}